\documentclass[12pt]{article}

\usepackage{graphicx}
\usepackage{bbm}
\usepackage{amsmath}
\usepackage{amsfonts}
\usepackage{amsthm}
\usepackage{amssymb}

\usepackage{color}

        \theoremstyle{plain}
        
        \newtheorem{proposition}{Proposition}[section]

        \theoremstyle{remark}
        \newtheorem{remark}{\bf Remark}[section]
        \theoremstyle{remark}
        
        \theoremstyle{remark}
        
\newcommand{\half}{\frac{1}{2}}

\newcommand{\demi}{\frac{1}{2}}

\newcommand{\R}{{\mathbb{R}}}

\newcommand{\dt}{\partial_t}
\newcommand{\dx}{\partial_x}

\newcommand{\eps}{\varepsilon}

\newcommand{\cint}[1]{\langle #1 \rangle}
\newcommand{\av}[1]{\cint{#1}}

\newcommand{\tn}{t_n}
\newcommand{\tnpun}{t_{n+1}}
\newcommand{\imax}{i_{max}}

\newcommand{\xipdemi}{x_{i+\demi}}
\newcommand{\ximdemi}{x_{i-\demi}}

\newcommand{\fni}{f_i^n}
\newcommand{\fnpuni}{f_i^{n+1}}
\newcommand{\fnimun}{f_{i-1}^n}
\newcommand{\fnipun}{f_{i+1}^n}

\newcommand{\One}{\mathbbm{1}}
\newcommand{\deltax}{\delta_x}

\newcommand{\sigmai}{\sigma_{i}}
\newcommand{\sigmapdemi}{\sigma_{i+\demi}}
\newcommand{\sigmamdemi}{\sigma_{i-\demi}}
\newcommand{\alphai}{\alpha_{i}}
\newcommand{\alphapdemi}{\alpha_{i+\demi}}

\newcommand{\phipdemi}{\phi_{i+\frac{1}{2}}}
\newcommand{\phimdemi}{\phi_{i-\frac{1}{2}}}
\newcommand{\Phipdemi}{\Phi_{i+\frac{1}{2}}}
\newcommand{\Phimdemi}{\Phi_{i-\frac{1}{2}}}
\newcommand{\Dx}{\Delta x}
\newcommand{\Dt}{\Delta t}

\renewcommand{\R}{{\mathbb{R}}}
\renewcommand{\eps}{\varepsilon}

\renewcommand{\dt}{\partial_t}

\renewcommand{\Dt}{\Delta t}
\renewcommand{\demi}{\frac{1}{2}}

\newcommand{\rhonpun}{\rho^{n+1}}
\newcommand{\rhon}{\rho^n}
\newcommand{\rhoni}{\rhon_i}
\newcommand{\rhonimun}{\rhon_{i-1}}
\newcommand{\rhonipun}{\rhon_{i+1}}
\newcommand{\rhonpuni}{\rhonpun_i}
\newcommand{\rhonpunimun}{\rhonpun_{i-1}}
\newcommand{\rhonpunipun}{\rhonpun_{i+1}}
\newcommand{\rhonipdemi}{\rhon_{i+\demi}}
\newcommand{\rhonLipdemi}{\rho^{nL}_{i+\demi}}
\newcommand{\rhonRipdemi}{\rho^{nR}_{i+\demi}}

\newcommand{\modif}{}

\author{Luc Mieussens}
\date{}                                           % Activate to display a given date or no date

\begin{document}
% \begin{flushright}
% {\it draft, for submission to the Journal of Computational Physics, \today}
% \end{flushright}

\begin{center}

{\bf On the Asymptotic Preserving property of the
  Unified Gas Kinetic Scheme for the diffusion limit of linear
  kinetic models \\}

\vspace{1cm}

Luc Mieussens\footnote{Univ. Bordeaux, IMB, UMR 5251, F-33400 Talence, France.\\
CNRS, IMB, UMR 5251, F-33400 Talence, France.\\
INRIA, F-33400 Talence, France.\\
({\bf  \tt Luc.Mieussens@math.u-bordeaux1.fr})}

\end{center}

\bigskip

{\bf Abstract. } 
The unified gas kinetic scheme (UGKS) of K. Xu et al.~\cite{XH2010},
originally developed for multiscale gas dynamics problems, is applied
in this paper to a linear kinetic model of radiative transfer
theory. While such problems exhibit purely diffusive behavior in the
optically thick (or small Knudsen) regime, we prove that UGKS is still
asymptotic preserving (AP) in this regime, but for the free
transport regime as well. Moreover, this scheme is modified to include
a time implicit discretization of the limit diffusion equation, and
to correctly capture the solution in case of boundary layers. Contrary
to many AP schemes, this method is based on a standard finite volume
approach, it does neither use any decomposition of the solution, nor
staggered grids. Several numerical tests demonstrate the properties of
the scheme. 

\bigskip


%--------------------------------------------------------------------------
\section{Introduction}

%--------------------------------------------------------------------------

Kinetic models are efficient tools to describe the dynamics of
systems of particles, like in rarefied gas dynamics (RGD), neutron
transport, semi-conductors, or radiative transfer. Numerical
simulations based on these models require important computational
resources, but modern computers make it possible to simulate
realistic problems.

These simulations can be made much faster when the ratio between the
mean free path of particles and a characteristic macroscopic length
(the so-called Knudsen number in RGD, denoted by $\eps$ in this paper)
is small. In such cases, the system of particles is accurately
described by a macroscopic model (Euler or Navier-Stokes equations in
RGD, diffusion equations in neutron or photon transport) that can be
numerically solved  much faster than with kinetic models. 

However, there are still important problems in which the numerical
simulation is difficult: in multiscale problems, $\eps$ can be very
small in some zones, and very large elsewhere (opaque vs. transparent
regions in radiative transfer). Standard numerical methods for kinetic
equations are very expensive in such cases, since, for stability and
accuracy reasons, they must resolve the smallest microscopic scale,
which is computationally expensive in small $\eps$ zones. By contrast,
macroscopic solvers are faster but may be inaccurate in large $\eps$ zones.

This is why multiscale numerical methods have been presented in the
past 20 years: the asymptotic-preserving (AP)
schemes. Such schemes are uniformly stable with respect to $\eps$
(thus their computational complexity does not depend on $\eps$), and
are consistent with the macroscopic model when $\eps$ goes to $0$ (the
limit of the scheme is a scheme for the macroscopic model). 

AP schemes have first been studied (for steady problems) in neutron
transport by Larsen, Morel and Miller~\cite{LMM1987}, Larsen and
Morel~\cite{LM1989}, and then by Jin and
Levermore~\cite{JL1991,JL1993}. For unstationary problems, the
difficulty is the time stiffness due to the collision operator. To
avoid the use of expensive fully implicit schemes, several
semi-implicit time discretizations schemes, based on a decomposition
of the distribution function between an equilibrium part and its
deviation, have been proposed by Klar~\cite{klar_sinum_1998}, and Jin,
Pareschi and Toscani~\cite{JPT2000} (see preliminary works
in~\cite{JPT1998,jin_sisc_1999} and extensions
in~\cite{JP2000b,JP2000,NP1998,klar_sinum_1999,klar_sisc_1999}). Similar
ideas have also been used by Buet et al. in~\cite{BCLM}, Klar and
Schmeiser~\cite{KS2001}, Lemou and Mieussens~\cite{LM2007,BLM2008},
and Carrillo et al.~\cite{CGL2008,CGLV2008}. The theory of well
balanced schemes is another way to obtain AP schemes, as in the work
of Gosse an Toscani~\cite{GT2002,GT2003}. Other approaches have been
recently proposed by Lafitte and Samaey~\cite{LS2012} and
Gosse~\cite{Gosse_2011}, but there extensions to more complex cases is
not clear. Finally, the idea of~\cite{jin_levermore} has been renewed
to obtain an AP scheme for linear equations on two-dimensional
unstructured meshes in the work of Buet, Despr\'es, and
Franck~\cite{BDF2012}. All these methods have advantages and
drawbacks, and there is still a need for other AP schemes.

A rather different approach has recently been proposed by K. Xu and
his collaborators, in the context of rarefied gas dynamics~\cite{XH2010}. This
method is called {unified gas kinetic scheme} (UGKS) and is based on
a gas kinetic scheme which has been developed by K. Xu since 2000
(see~\cite{Xu_2001} for the first reference and many other references
in~\cite{XH2010}). Roughly speaking, the UGKS is based on a finite volume
approach in which the numerical fluxes contain information from the
collision operator. In some sense, it has some connexions with the
well balanced schemes developed for hyperbolic problems with source
terms in~\cite{jin_levermore,GT2002,GT2003}, even if the
construction is completely different. While this approach to design AP
schemes looks very promising, it seems that it has not yet received the
attention it deserves from the kinetic community. This is probably due
to the fact that the nice properties of the UGKS presented
in~\cite{XH2010} are difficult to understand for people who are not
specialist of gas kinetic schemes. 

However, we believe that the UGKS approach is very general and can
benefit to many different kinetic problems. Let us mention that the
big advantage of the UGKS with respect to other methods is that is does
not require any decomposition of the distribution function (hence
there is no problem of approximation of the boundary conditions), it
does not use staggered grids (which is simpler for multi-dimensional
problems), and it is a finite volume method (there is no need of
discontinuous Galerkin schemes that are more expensive).

In this paper, our first goal is to present UGKS in a very simple
framework, so that it can be understood by any researcher interested
in numerical method for kinetic equations. We also want to show that
the UGKS can be successfully applied to other fields than RGD. Here, it is
used to design an AP scheme for linear kinetic equations, namely a
simple model of radiative transfer. Such an extension is not obvious,
since linear models exhibit a purely diffusive (parabolic) behavior in
the small $\eps$ regimes, while models from RGD (like the Boltzmann
equations) have a rather convection (hyperbolic) behavior. Indeed,
even if the UGKS is originally made to correctly describe this convection regime
and to capture the small viscous effects (like in the compressible
Navier-Stokes equations), we prove in this paper that it can also capture a
purely diffusive effect. Moreover, we propose
several extensions: implicit diffusion, correct boundary conditions
for boundary layers, treatment of collision operator with non
isotropic scattering kernel. The scheme is proved to be AP in both
free transport and diffusion regimes, and is validated with several
numerical tests.

The outline of our paper is the following. In section~\ref{sec:ugks},
we present the linear kinetic model, and its approximation by the
UGKS. Its asymptotic properties are analyzed in
section~\ref{sec:analysis}. Some extensions are given in
section~\ref{sec:extensions}, and the scheme is validated with various
numerical tests in section~\ref{sec:res}.

%--------------------------------------------------------------------------
\section{The UGKS for a linear transport model}
\label{sec:ugks}
%--------------------------------------------------------------------------

%=========
   \subsection{A linear transport model and its diffusion limit}
   \label{subsec:}
%=========

The linear transport equation is a model for the evolution (by
transport and interaction) of particles in
some medium. In this paper, we are mainly concerned by the radiative
transfer equation, which reads
{\modif
\begin{equation*}
\frac{1}{c}\dt \phi + \Omega\cdot  \nabla_{r} \phi =
\sigma(\frac{1}{4\pi}\int \phi \,d\Omega - \phi) -\alpha  \phi + G,  
\end{equation*}
where $\phi(t,r,\Omega)$ is the spectral intensity in the
position-direction phase space that depends on time $t$, position
$r=(x,y,z)\in \R^3$, and angular direction of propagation of particles
$\Omega\in S^2$, while $c$ is their velocity (the speed of
light). Moreover, $\sigma$ is the scattering cross section, $\alpha$
is the absorption cross section, and $G$ is an internal source of
particles. These three last quantities may depend on $x$, but they are
independent of $\Omega$. The linear operator $\phi\mapsto
\frac{1}{4\pi}\int \phi\,d\Omega - \phi$ models the scattering of the
particles by the medium and acts only on the angular dependence of
$\phi$. This simple model does not allow for particles of possibly
different energy (or frequency); it is called ``one-group'' or
``monoenergetic'' equation.  

In order to study the diffusion regime
corresponding to this equation, a standard dimensional analysis is
made (see~\cite{LMM1987} for details). We choose a macroscopic length
scale $r^*$, like the size of the computational domain. We assume that
this length is much larger than the typical mean free path $\lambda^*$
(defined by a typical value of $1/\sigma$), and we denote by
$\eps$ the ratio $\lambda^*/r^*$ which is supposed to be much smaller
than $1$. We choose a macroscopic time scale $t^*$
which is much larger than the typical mean free time
$\tau^*=\lambda^*/c$, so that $\frac{\tau^*}{t^*}=\eps^2$. Finally, we
assume that the absorption cross section $\alpha$ and the source $G$
are of the order $O(\eps^2)$ as compared to $\sigma$. Then, with the
non dimensional scaled variables $t'=t/t^*$, $r'=r/r^*$,
$\sigma'=\sigma/\sigma^*$, $\alpha'=\alpha/\alpha^*$, $G'=G/G^*$, we
get the following equation
\begin{equation*}
\eps \partial_{t'} \phi + \Omega\cdot  \nabla_{r'} \phi =
\frac{\sigma'}{\eps}(\frac{1}{4\pi}\int \phi \,d\Omega - \phi) -\eps\alpha'  \phi + \eps G'.
\end{equation*}
In the following, we drop all the ' in the equations, since we always
work with the non dimensional variables.}

In this paper, we consider this one-group equation in the slab
geometry: we assume that $\phi$  depends only on the slab axis
variable $x\in \R$. Then it can be shown that the average of $\phi$
with respect to the $(y,z)$ cosine directions of $\Omega$, denoted by
$f(t,x,v)$, satisfies the one-dimensional equation
\begin{equation}  \label{eq-onegroup}
\eps \dt f + v \dx f = \frac{\sigma}{\eps}(\av{f}-f) -\eps\alpha  f + \eps
G,
\end{equation}
where $v\in[-1,1]$ is the $x$ cosine direction of $\Omega$ and the
operator $\av{.}$ is such that $\av{\phi}=\half\int_{-1}^1\phi(v)\,
dv$ is the average of every $v$-dependent function $\phi$.

% We assume that the
% cross sections satisfy the inequalities $0<\sigma_m\leq \sigma(x)\leq
% \sigma_M$ and $0\leq \sigma_A(x)\leq \sigma_{AM}$ for every $x$. 

When $\eps$ becomes small, it is well known that the solution $f$
of~(\ref{eq-onegroup}) tends to its own
average density $\rho=\av{f}$, which is a solution of the asymptotic
diffusion limit
\begin{equation}  \label{eq-diff}
\dt \rho -\dx \kappa  \dx\rho =-\alpha\rho +G,
\end{equation}
where the diffusion coefficient is
$\kappa(x)=\frac{\av{v^2}}{\sigma(x)}=\frac{1}{3\sigma(x)}$. An asymptotic preserving
scheme for the linear kinetic equation~(\ref{eq-onegroup}) is a
numerical scheme that discretizes~(\ref{eq-onegroup}) in such a way
that it leads to a correct discretization of the diffusion
limit~(\ref{eq-diff}) when $\eps$ is small.

%=========
\subsection{First ingredient of the UGKS: a finite volume scheme}
\label{subsec:FV}
%=========

The first ingredient of the UGKS is a finite volume
approach. Equation~(\ref{eq-onegroup}) is integrated over a time
interval $[t_n,t_{n+1}]$ and over a space cell $[\xipdemi,\ximdemi]$ to
obtain the following relation:
\begin{equation*}
\frac{\fnpuni-\fni}{\Dt}  + \frac{1}{\Dx}\left(
  \phipdemi - \phimdemi
  \right) 
 = \frac{1}{\Dt\Dx}\int_{t_n}^{\tnpun} \int_{\ximdemi}^{\xipdemi} (\frac{\sigma}{\eps^2}(\rho-f)-\alpha f )   \, dxdt + G,
\end{equation*}
where $f^n_i=\frac{1}{\Dx}\int_{x_{i-\demi}}^{x_{i+\demi}}
f(t_n,x,v)\, dv$ is the average of $f$ over a space cell and
$\phipdemi$ is the microscopic flux across the interface $\xipdemi$:
\begin{equation}\label{eq-flux_cont} 
  \phipdemi= \frac{1}{\eps\Dt}\int_{t_n}^{\tnpun} v f(t,\xipdemi,v)  \, dt. 
\end{equation}
To obtain a scheme which is uniformly stable with respect to $\eps$, the collision term, which is the stiffest term in the previous relation when $\eps$ is small, must be discretized by an implicit approximation. For simplicity, we use here the standard right-rectangle quadrature which gives a first order in time approximation:
\begin{equation*}
  \frac{1}{\Dt\Dx}\int_{t_n}^{\tnpun} \int_{\ximdemi}^{\xipdemi} (\frac{\sigma}{\eps^2}(\rho-f)-\alpha f )   \, dxdt 
\approx \frac{\sigmai}{\eps^2}(\rhonpuni-\fnpuni)-\alphai \fnpuni. 
\end{equation*}
We also assume that the scattering and absorption coefficient do not vary to much inside a cell, so that the average of the products $\sigma f$ and $\alpha f$ are close to the product of the averages of each terms. Taking the absorption term implicit is not necessary, and an explicit approximation could also be used.

Then the finite volume scheme for~(\ref{eq-onegroup}) reads:
\begin{equation}\label{eq-fv}
\frac{\fnpuni-\fni}{\Dt}  + \frac{1}{\Dx}\left(
  \phipdemi - \phimdemi
  \right) 
 = \frac{\sigmai}{\eps^2}(\rhonpuni-\fnpuni)-\alphai \fnpuni + G.
\end{equation}
Now, it remains to approximate the flux $\phipdemi$, and hence the value of $f$ at any time $t$ between $t_n$ and $t_{n+1}$ at the cell interface $\xipdemi$, by using the averaged values $\fni$, $\fnipun$, $\fnimun$, etc. The choice made for this approximation is the core of the UGKS scheme.

Note that in this construction, $v$ is not discretized: indeed, while
it can be discretized by any method, this would not change our
analysis. Keeping $v$ continuous at this stage considerably simplifies
the notations, what we do here. See section~\ref{sec:res} for an
example of velocity discretization.

%=========
\subsection{Second ingredient: a characteristic based approximation of the cell interface value $f(t,\xipdemi,v)$}
%=========

In this section, we explain how $f(t,x_{i+\demi},v)$ is reconstructed
in the flux defined in~(\ref{eq-flux_cont}). To emphasize the importance of
this approximation, note that the standard first-order upwind
approximation
\begin{equation*}
  f(t,\xipdemi,v) = \fni\One_{v>0} + \fnipun\One_{v<0},
\end{equation*}
where $\One_{v\lessgtr 0}=1$ if $v\lessgtr 0$ and $0$ else, is not a
choice that gives an AP scheme. Indeed, it introduces numerical
dissipation that dominates the physical diffusion in the diffusive
limit regime. This can also be interpreted as
follows: this approximation is nothing but the solution at time $t$ of the Riemann problem 
\begin{equation*}
  \begin{split}
& \dt f + \frac{1}{\eps}v\dx f = 0,    \\
& f(t_n,x,v)=\fni \text{ if } x>\xipdemi, \fnipun \text{ else}.  
\end{split}
\end{equation*}
This problem does not take into account the right-hand side of the real problem~(\ref{eq-onegroup}). 
It is well known that hyperbolic problem with source term cannot be
accurately approximated if the numerical flux is constructed by
ignoring the source term, in particular in limit regimes
(see~\cite{jin_levermore}, for instance). While it is not always easy
to take the source term into account in the numerical flux, the transport+relaxation structure of~(\ref{eq-onegroup}) makes this problem particularly simple.

Indeed, in~\cite{XH2010}, Xu and Huang propose to use the integral representation of the solution of the BGK equation, which, for our model~(\ref{eq-onegroup}), is obtained as follows. If the coefficients $\sigma$ and $\alpha$ are constant in space and time,
equation~(\ref{eq-onegroup}) is equivalent to 
\begin{equation*}
\frac{d}{dt}e^{\nu t}f(t,x+\frac{v}{\eps}t,v)
=e^{\nu  t}\Bigl(\frac{\sigma}{\eps^2}\rho(t,x+\frac{v}{\eps}t)+G\Bigr)  ,
\end{equation*}
where $\nu=\frac{\sigma}{\eps^2}+\alpha$. If they are not constant but
slowly varying in one cell, we consider this relation as an
approximation around $t_n$ and $\xipdemi$, and we denote by
$\sigma_{i+\demi}$, $\alpha_{i+\demi}$ and $\nu_{i+\demi}$ the
corresponding constant values of $\sigma$, $\alpha$ and $\nu$. Then we integrate this relation between $t_n$ and some $t<t_{n+1}$, and we replace $x+\frac{v}{\eps}t$ by $\xipdemi$ to get the following relation: 
\begin{equation}\label{eq-ftxipdemi}
\begin{split}
  f(t,\xipdemi,v) & \approx   e^{-\nu_{i+\demi}(t-\tn)}f(\tn,\xipdemi-\frac{v}{\eps}(t-t_n)) \\
& \quad + \int_{t_n}^{t} e^{-\nu_{i+\demi}(t-\tn)}\frac{\sigma_{i+\demi}}{\eps^2}\rho(s,\xipdemi-\frac{v}{\eps}(t-s))  \, ds \\
& \quad +\frac{1-e^{-\nu_{i+\demi}(t-\tn)}}{\nu_{i+\demi}} G.
\end{split}
\end{equation}

Now, it remains to design an approximation of two terms: the first one
is $f$ at $t_n$ around $\xipdemi$, that is to say
$f(\tn,\xipdemi-\frac{v}{\eps}(t-t_n)) $, and the second one is $\rho$
between $t_n$ and $t_{n+1}$ around $\xipdemi$, that is to say $\rho(s,\xipdemi-\frac{v}{\eps}(t-s)) $.

To approximate $f$ at $t_n$ around $\xipdemi$, the simplest approach
is to use a piecewise constant reconstruction: 
\begin{equation}  \label{eq-piecewisef}
f(t_n,x,v)=\left\lbrace
\begin{split}
& \fni \quad \text{ if } \quad x<\xipdemi \\
& \fnipun \quad \text{ if } \quad x>\xipdemi.
\end{split}
\right.
\end{equation}
Of course, a more accurate reconstruction can be obtained (for
instance a piecewise linear reconstruction with slope limiters, as
in~\cite{XH2010}), but the presentation of the scheme is simpler with
this zeroth order reconstruction.

This is for the approximation of $\rho$ between $t_n$ and $t_{n+1}$
around $\xipdemi$ that we need the second main idea of K. Xu: this
reconstruction is piecewise {\it continuous}. This can be surprising,
since $\rho$ is the velocity average of $f$ which is represented by
piecewise discontinuous function, but this is the key idea that allows
the scheme to capture the correct diffusion terms in the small $\eps$
limit. In~\cite{XH2010,HXY2012}, this reconstruction is piecewise linear in
space and time. However, in the context of the diffusion limit, we
found that the a piecewise constant in time reconstruction is
sufficient, while a piecewise continuous linear reconstruction in space
is necessary. First, we define the unique interface value $\rhonipdemi$
at time $t_n$ by:
\begin{equation}  \label{eq-rhonipdemi}
\rhonipdemi= \cint{\fni\One_{v>0} + \fnipun\One_{v<0}}.
\end{equation}
Then, we define the following reconstruction for $t$ in
$[t_n,t_{n+1}]$ and $x$ around $\xipdemi$:
\begin{equation}\label{eq-rhoreconst} 
\rho(t,x)=\left\lbrace
\begin{split}
& \rhonipdemi + \deltax\rhonLipdemi(x-\xipdemi) \quad \text{ if } x<\xipdemi \\
& \rhonipdemi + \deltax\rhonRipdemi(x-\xipdemi) \quad \text{ if } x>\xipdemi \\
\end{split}
\right.
\end{equation}
with left and right one-sided finite differences: 
\begin{equation}  \label{eq-slopes}
\deltax\rhonLipdemi=\frac{\rhonipdemi-\rhoni}{\Dx/2} \quad \text{ and
} \quad 
\deltax\rhonRipdemi=\frac{\rhonipun-\rhonipdemi}{\Dx/2}
\end{equation}

\begin{remark}
  When $\eps$ is very small, the foot of the characteristic
  $\xipdemi-\frac{v}{\eps}(t-s)$ in~(\ref{eq-ftxipdemi}) can be very
  far from $\xipdemi$, and hence using the
  reconstruction~(\ref{eq-rhoreconst}) might be very
  inaccurate. However, note that in~(\ref{eq-ftxipdemi}), $\rho$ is
  multiplied by an exponential term which is very small in this case,
  and we can hope that the inaccuracy made in the
  reconstruction~(\ref{eq-rhoreconst}) has not a too strong
  influence. Indeed, our asymptotic analysis and our numerical tests
  show that the accuracy of the scheme is excellent, see
  sections~\ref{sec:analysis} and~\ref{sec:res}.
\end{remark}

%=========
   \subsection{Numerical flux}
   \label{subsec:numflux}
%=========

Now, the numerical flux $\phipdemi=
\frac{1}{\eps\Dt}\int_{t_n}^{\tnpun} v f(t,\xipdemi,v)  \, dt $ can be
computed exactly by using expressions~(\ref{eq-ftxipdemi},\ref{eq-piecewisef},\ref{eq-rhoreconst})
to get
\begin{equation}\label{eq-flux}
\begin{split}
 \phipdemi& = A_{i+\demi}v\left( \fni \One_{v>0}
                                 +  \fnipun\One_{v<0}\right)     \\
& \quad + C_{i+\demi}v\rhonipdemi \\
& \quad
 + D_{i+\demi}v^2( \deltax\rhonLipdemi\One_{v>0} + \deltax\rhonRipdemi\One_{v<0})\\
& \quad + E_{i+\demi}vG
\end{split}
\end{equation}
where the coefficients $A_{i+\demi}=A(\Dt, \eps, \sigmapdemi,\alphapdemi)$, $C_{i+\demi}=C(\Dt, \eps, \sigmapdemi,\alphapdemi)$, $D_{i+\demi}=D(\Dt, \eps, \sigmapdemi,\alphapdemi)$, and
$E_{i+\demi}=E(\Dt, \eps, \sigmapdemi,\alphapdemi)$ are defined by the
following functions: 
\begin{align}
  A(\Dt, \eps, \sigma,\alpha)& =\frac{1}{\Dt \eps \nu}(1-e^{-\nu \Dt}) \\
  C(\Dt, \eps, \sigma,\alpha)& =\frac{\sigma}{\Dt \eps^3 \nu}
                                \Bigl(\Dt-\frac{1}{\nu}\bigl(1-e^{-\nu \Dt}\bigr)\Bigr)
  \\
  D(\Dt, \eps, \sigma,\alpha)& =-\frac{\sigma}{\Dt \eps^4 \nu^2}
            \Bigl(\Dt\bigl(1+e^{-\nu \Dt}\bigr)-\frac{2}{\nu}\bigl(1-e^{-\nu \Dt}\bigr)\Bigr)
  \\
  E(\Dt, \eps, \sigma,\alpha)& =\frac{1}{\Dt \eps \nu}
                                \Bigl(\Dt-\frac{1}{\nu}\bigl(1-e^{-\nu
                                  \Dt}\bigr)\Bigr)\label{eq-coeff_flux}
\end{align}
where we remind that $\nu=\frac{\sigma}{\eps^2}+\alpha$.

Now, using~(\ref{eq-fv}), $\fnpuni$ can be obtained, provided that we
can first determine $\rhonpuni$. This is done by using the
conservation law in the following section.

%=========
   \subsection{Conservation law}
%   \label{subsec:}
%=========
It is now well known that semi-implicit schemes for relaxation kinetic
equations can be solved explicitly by using the corresponding
discrete conservation laws: see~\cite{PP2007} where it was first
suggested for the BGK equation, and~\cite{XH2010} for the use of this
technique to design the UGKS.

The idea is to eliminate $\fnpuni$ in~(\ref{eq-fv}) by taking its
$v$-average, hence obtaining the following discrete conservation law
\begin{equation}\label{eq-conslaws}
  \frac{\rhonpuni-\rhoni}{\Dt}  + \frac{1}{\Dx}\left(
  \Phipdemi - \Phimdemi
  \right) 
 = - \alpha_i \rhonpuni +  G
\end{equation}
where the macroscopic numerical flux is
\begin{equation*}
  \Phipdemi = \cint{\phipdemi}.
\end{equation*}
By using~(\ref{eq-flux}), we find: 
\begin{equation}\label{eq-macrosflux}
 \Phipdemi = A_{i+\demi}\cint{v\fni\One_{v>0}+ v\fnipun\One_{v<0}}      
+ D_{i+\demi}\frac{1}{3}\frac{\rhonipun-\rhoni}{\Dx}.
\end{equation}
Note that relation~(\ref{eq-conslaws}) is still implicit, but it can
be solved explicitly.

%=========
   \subsection{Summary of the numerical scheme}
   \label{subsec:summary}
%=========

Finally, $\fnpuni$ is computed as follows for every cell $i$: 
  \begin{enumerate}

  \item compute $\rhonpuni$ by solving:
\begin{equation}\label{eq-UGKS_rho} 
  \frac{\rhonpuni-\rhoni}{\Dt}  + \frac{1}{\Dx}\left(
  \Phipdemi - \Phimdemi
  \right) 
 = - \alpha_i \rhonpuni +  G,
\end{equation}

  \item compute $\fnpuni$ by solving: 
\begin{equation}\label{eq-UGKS_f} 
\frac{\fnpuni-\fni}{\Dt}  + \frac{1}{\Dx}\left(
  \phipdemi - \phimdemi
  \right) 
 = \frac{\sigma_i}{\eps^2}(\rhonpuni-\fnpuni)- \alpha_i \fnpuni + G,
\end{equation}
  \end{enumerate}
where macroscopic and microscopic fluxes are given
by~(\ref{eq-macrosflux}) and~(\ref{eq-flux}).

Of course, this scheme must be supplemented by numerical boundary
conditions. This will be detailed in the next section, after the
asymptotic analysis.

%--------------------------------------------------------------------------
   \section{Asymptotic analysis}
%--------------------------------------------------------------------------
\label{sec:analysis}

%=========
   \subsection{Free transport regime}
%   \label{subsec:}
%=========
The behavior of the scheme in the small $\sigma$, $\alpha$ limit is completely
determined by the following property of the coefficient functions $A$,
$C$, $D$, and $E$.
\begin{proposition}\label{prop:limit_sigma_coef} 
  When $\sigma$ and $\alpha$ tend to 0 (while $\eps$ is constant), we have: 
\begin{itemize}
  \item  $A(\Dt,\eps,\sigma,\alpha)$ tends to $\frac{1}{\eps}$
  \item $C$ and $D(\Dt,\eps,\sigma,\alpha)$ tend to 0
  \item $E(\Dt,\eps,\sigma,\alpha)$ tends to $\frac{\Dt}{2}$.
\end{itemize}
\end{proposition}
As a consequence, the microscopic flux $\phipdemi$ defined
in~(\ref{eq-flux}) has the following limit:
\begin{equation*}
  \phipdemi \xrightarrow[\sigma\rightarrow 0]{}
  \frac{v}{\eps}(\fni\One_{v>0}+\fnipun\One_{v<0}) + \frac{\Dt}{2}vG.
\end{equation*}
This is nothing but a consistent first-order upwind flux (plus a
constant term that has no influence in the scheme), and the limit finite
volume scheme~(\ref{eq-fv}) is
\begin{equation*}
  \frac{\fnpuni-\fni}{\Dt}  + \frac{1}{\Dx}\frac{v}{\eps}
\left(
(\fni\One_{v>0}+\fnipun\One_{v<0}) 
- 
(\fnimun\One_{v>0}+\fni\One_{v<0})
  \right) 
 = G,
\end{equation*}
which is indeed a consistent approximation of the limit transport
equation~(\ref{eq-onegroup}) when $\sigma$ and $\alpha$ tend to 0. As a
consequence, the UGKS is AP in this limit.

%=========
   \subsection{Diffusion regime}
   \label{subsec:diffusion_regime}
%=========
Similarly, the behavior of the scheme in the small $\eps$ limit is completely
determined by another property of the coefficient functions $A$
and $D$. 
\begin{proposition}\label{prop:limit_eps_coef} 
  When $\eps$ tends to 0, we have: 
\begin{itemize}
  \item  $A(\Dt,\eps,\sigma,\alpha)$ tends to 0
%  \item $C(\Dt,\eps,\sigma,\alpha)$ is asymptotically equivalent to $\frac{1}{\eps}$
  \item $D(\Dt,\eps,\sigma,\alpha)$ tends to $-\frac{1}{\sigma}$.
\end{itemize}
\end{proposition}
As a consequence, the macroscopic flux $\Phipdemi$ defined
in~(\ref{eq-macrosflux}) has the following limit:
\begin{equation}\label{eq-limitmacrosflux} 
  \Phipdemi \xrightarrow[\eps\rightarrow 0]{}
  -\frac{1}{3\sigmapdemi}\frac{\rhonipun-\rhoni}{\Dx},
\end{equation}
and hence the limit of the discrete conservation
law~(\ref{eq-conslaws}) is
\begin{equation}\label{eq-num-diffusion-limit}
  \frac{\rhonpuni-\rhoni}{\Dt}   -  \frac{1}{\Dx}
\biggl(
 \frac{1}{3\sigmapdemi}\frac{\rhonipun-\rhoni}{\Dx}
- 
\frac{1}{3\sigmamdemi}\frac{\rhoni-\rhonimun}{\Dx}
 \biggr) =- \alpha_i \rhonpuni + G
\end{equation}
which is a consistent approximation of the diffusion
limit~(\ref{eq-diff}), with a standard three-point centered
approximation of the second order derivative of $\rho$. This proves
that UGKS is AP for the limit $\eps \rightarrow 0$.

{\modif Moreover, it is interesting to look at the case where
  the scattering coefficients are discontinuous. In the derivation of
  the UGKS, the assumption that the scattering coefficients do not
  vary too much in and around a given cell is used several times: in
  the rectangle formula to obtain the right-hand side
  of~(\ref{eq-fv}), and in the characteristic based solution of the
  Riemann problem~(\ref{eq-ftxipdemi}). In case of discontinuous
  coefficients, the previous assumption is not satisfied, and the
  derivation is not correct. However, one can still use the same
  scheme, and see how the interface values $\sigma_{i+\demi}$ and
  $\alpha_{i+\demi}$ can be defined to obtain correct results.

  Assume that $\sigma$ and $\alpha$ are piecewise continuous, that is
  to say that they are continuous in each cell (with cell averages
  denoted by $\sigma_i$ and $\alpha_i$), and have possible
  discontinuities across cell interfaces. If we simply define
  $\sigma_{i+\demi}$ as the arithmetic average of $\sigma$ in the two
  adjacent cells, that is to say
  $\sigma_{i+\demi}=(\sigma_i+\sigma_{i+1})/2$ (and the same for
  $\alpha$), then scheme~(\ref{eq-UGKS_rho})--(\ref{eq-UGKS_f}) is
  unchanged, and its AP property is still satisfied.

Let us have a look to the asymptotic scheme we get in the diffusion
limit: as found above, the scheme converges to~(\ref{eq-num-diffusion-limit}).
It means that the diffusion coefficient--which is
$\kappa(x)=1/(3\sigma(x))$ in the continuous case--is approximated at
the cell interface by
$\kappa_{i+\demi}=1/(3\sigma_{i+\demi})$. If the definition of
$\sigma_{i+\demi}$ is injected into this relation, it is easy to see that $\kappa_{i+\demi}$ can be written as
\begin{equation*}
  \kappa_{i+\demi}=2\frac{\kappa_i\kappa_{i+1}}{\kappa_i+\kappa_{i+1}},
\end{equation*}
which is the harmonic average of $\kappa$ over the two adjacent
cells. Consequently, we recover a standard finite volume scheme for
the diffusion equation with discontinuous coefficients, which is known
to be second order accurate (see~\cite{wesseling}, for instance).
}

%=========
   \subsection{Boundary conditions}
   \label{subsec:BC}
%=========

If we consider equation~(\ref{eq-onegroup}) for $x$ in the bounded
domain $[0,1]$, we need the following boundary conditions:
\begin{equation}  \label{eq-BC}
f(t,x=0,v>0)=f_L(v) \quad \text{ and } \quad f(t,x=1,v<0)=f_R(v),
\end{equation}
where $f_L$ and $f_R$ can be used to model inflow or reflexion
boundary conditions. In the case of an inflow boundary condition, if
$f_L$ and $f_R$ are independent of $v$, the diffusion limit is
still~(\ref{eq-diff}) with the corresponding Dirichlet boundary data
$\rho(t,0)=f_L$ and $\rho(t,1)=f_R$. However, when one of the boundary
data (say $f_L$) depends on $v$ (which is called a non isotropic
data), using this diffusion approximation requires some
modifications. First, a boundary layer corrector must be added to
$\rho$ to correctly approximate $f$. Moreover, $f$ is well
approximated by the solution $\rho$ of the diffusion equation outside
the boundary layer only if this equation is supplemented by the boundary
condition 
\begin{equation}  \label{eq-CLchandra}
\rho(t,x=0)=\int_0^1W(v)f_L(v) \, dv= 2\cint{W f_L \One_{v>0}}
\end{equation}
where $W(v)$ is a special function that can be well approximated by
$0.956 v + 1.565 v^2\approx 3/2v^2+v$
(see~\cite{LMM1987}).

Now, we explain how these boundary conditions can be taken into
account numerically in the UGKS. We assume that we have $\imax$ space
cells in $[0,1]$. As compared to our derivation of
scheme~(\ref{eq-fv}), the difference is that the derivation of the
numerical fluxes at the boundaries $\phi^n_{\demi}$ and
$\phi^n_{\imax+\demi}$ must take into account the boundary data. For
simplicity, we only describe in details how the left boundary flux is
constructed. The definition~(\ref{eq-ftxipdemi}) is still valid, but the
integral representation of $f$ at this left boundary now is
\begin{equation*}
  f(t,x_{\demi},v)  =\left\lbrace
    \begin{split}
& f_L \qquad \text{if } v>0  \\
& e^{-\nu_{\demi}(t-\tn)}f(\tn,x_{\demi}-\frac{v}{\eps}(t-t_n)) 
%\\
%& \qquad 
+ \int_{t_n}^{t} e^{-\nu_{\demi}(t-\tn)}\frac{\sigma}{\eps^2}\rho(s,x_{\demi}-\frac{v}{\eps}(t-s))  \, ds \\
& \quad +\frac{1-e^{-\nu_{\demi}(t-\tn)}}{\nu_{\demi}} G \qquad \text{if } v<0      
    \end{split}
\right. ,
\end{equation*}
where $f(\tn,x_{\demi}-\frac{v}{\eps}(t-t_n))$ and $\rho(s,x_{\demi}-\frac{v}{\eps}(t-s))$
have to be defined for $v<0$ only, that is to say on the right-hand
side of the left boundary. According to the
approximations~(\ref{eq-piecewisef}) and~(\ref{eq-rhoreconst}), we set
for $v<0$ and $x>x_{\demi}$:
\begin{equation*}
  f(t_n,x,v)=f^n_1,
\end{equation*}
and for $t$ in $[t_n,t_{n+1}]$ 
\begin{equation*}
  \rho(t,x)=\rhon_{\demi}+\delta_x \rho^{nR}_{\demi}(x-x_{\demi}),
\end{equation*}
where according to~(\ref{eq-rhonipdemi}) the left boundary value
of $\rho$ now should be
\begin{equation}\label{eq-rhondemi} 
  \rho^n_{\demi}  = \cint{f_L\One_{v>0} + f^n_1\One_{v<0}},
\end{equation} 
and $\deltax\rho^{nR}_{\demi}=\frac{\rho^n_1-\rhon_{\demi}}{\Dx/2}$.
With this definitions, the numerical flux at the left boundary is found
to be
\begin{equation}\label{eq-fluxleft}
%\begin{split}
 \phi_{\demi} = \frac{v}{\eps}f_L \One_{v>0} 
+ A_{\demi}v f^n_1 \One_{v<0}
+ C_{\demi}v\rho^n_{\demi} \One_{v<0} 
+ D_{\demi}v^2 \deltax\rho^{nR}_{\demi}\One_{v<0}
+ E_{\demi}vG\One_{v<0},
%\end{split}
\end{equation}
where the coefficients $A_{\demi},B_{\demi},C_{\demi},D_{\demi}$ and
$E_{\demi}$ are defined as in section~\ref{subsec:numflux}. The
corresponding macroscopic flux is
\begin{equation}  \label{eq-macrosfluxleft}
%\begin{split}
 \Phi_{\demi} = \frac{1}{\eps}\cint{vf_L \One_{v>0}} 
+ A_{\demi}\cint{v f^n_1 \One_{v<0}}
+ \cint{v\One_{v<0}}C_{\demi}\rho^n_{\demi} 
+ \cint{v^2\One_{v<0}}D_{\demi} \deltax\rho^{nR}_{\demi}\One_{v<0}
+\cint{v\One_{v<0}} E_{\demi}G\One_{v<0}.
%\end{split}
\end{equation}
We could have used the exact values $\cint{v\One_{v<0}}=-1/4$ and
$\cint{v^2\One_{v<0}}=1/6$ in~(\ref{eq-macrosfluxleft}), but in
practice, when $v$ is discretized, it
is better to use the numerical approximation of $\cint{v\One_{v<0}}$ and
$\cint{v^2\One_{v<0}}$ corresponding to the chosen quadrature formula: this avoids an important loss of
accuracy when $\eps$ goes to 0.

With these relations that define $\phi_{\demi} $ and $\Phi_{\demi} $,
the numerical scheme summarized in section~\ref{subsec:summary} can be
used in every cell $i$ from $1$ to $i_{max}$. 

Unfortunately, this scheme is not uniformly stable for small $\eps$,
hence cannot be AP in the diffusion limit. The reason is that the
macroscopic flux now contains two unbounded terms that are
$\frac{1}{\eps}\cint{vf_L \One_{v>0}} $ and $-
\frac{1}{4}C_{\demi}\rho^n_{\demi} $ (we can prove that $C$ is
asymptotically equivalent to $\frac{1}{\eps}$). This fact is not
observed in~\cite{XH2010}, since UGKS is not applied in a diffusive
scaling, and hence the boundary data does not contribute as a $1/\eps$
term in the boundary flux. However, this drawback can be easily
corrected. First, we propose the following simple correction: it is
sufficient to modify the definition~(\ref{eq-rhondemi}) of the
boundary value $\rho^n_{\demi}$ to
\begin{equation}\label{eq-rhondemi_modified}
    \rho^n_{\demi}  = -\frac{\cint{vf_L\One_{v>0}}}{\cint{v\One_{v<0}}},
\end{equation}
so that the two unbounded terms exactly cancel each other. This
definition ensures that $\Phi_{\demi}$ is uniformly bounded with
respect to $\eps$, and it is reasonable to conjecture that the
corresponding scheme is uniformly stable. 

Now, we investigate the two numerical limits of our scheme with this
modified boundary condition. Note that this definition, while not
consistent with the physical boundary data,
is a standard approximation of the exact Dirichlet boundary condition
of the diffusion limit. This means that we can hope for a correct
behavior of the scheme in the diffusion limit. Indeed, following the
same analysis as in section~\ref{subsec:diffusion_regime}, we find
that the discrete conservation law~(\ref{eq-conslaws}) in the first
cell $i=1$ tends to 
\begin{equation}\label{eq-diff1} 
  \frac{\rho^{n+1}_1-\rho^n_1}{\Dt}   -  \frac{1}{\Dx}
\biggl(
 \frac{1}{3\sigma_{\frac{3}{2}}}\frac{\rho^n_2-\rho^n_1}{\Dx}
- 
\frac{1}{3\sigma_{\demi}}\frac{\rhon_1-\rhon_{\demi}}{\Dx}
 \biggr) =- \alpha_i \rho^{n+1}_1 + G,
\end{equation}
which is consistent approximation of the diffusion limit equation,
with a Dirichlet boundary data $\rho^n_{\demi}$ given
by~(\ref{eq-rhondemi_modified}). This boundary data is the correct one
if $f_L$ is isotropic (since we get $\rhon_{\demi}=f_L$), but is only
an (standard) approximation of the exact value~(\ref{eq-CLchandra}) if $f_L$ is
not isotropic. Note that a similar problem occurs with several AP
schemes in the literature, like in the schemes
of~\cite{JPT1998,LM2007}, while it seems to be greatly reduced in the
recent scheme of~\cite{LM2013}.  We propose in
section~\ref{subsec:ext_CL} a modification of the boundary condition
of our scheme to increase the accuracy near the boundaries.

However, contrary to the schemes that have been mentioned, the
inaccuracy of our modified boundary value $\rhon_{\demi}$ has no
influence in the free transport regime, which is another interesting
property of the UGKS. Indeed, since coefficients $C,D$, and $E$ tend
to 0 for small $\sigma$, the boundary fluxes $\phi_{\demi}$ and
$\Phi_{\demi}$ tend to the
simple first-order upwind boundary fluxes
\begin{align*}
&   \phi_{\demi}=\frac{v}{\eps}(f_L\One_{v>0} + f^n_1\One_{v<0}),  \\
& \Phi_{\demi}=\cint{\phi_{\demi}}=\frac{1}{\eps}\cint{vf_L\One_{v>0} + vf^n_1\One_{v<0}}.
\end{align*}

{\modif
%=========
   \subsection{Stability}
   \label{subsec:stab}
%=========
   Of course, the AP property also requires that the scheme is
   uniformly stable with respect to $\eps$, $\sigma$, and
   $\alpha$. Such a property is not easily derived, but looking at the
   asymptotic limits can give some indications. If the scheme is
   indeed uniformly stable with respect to $\eps$, then the CFL
   stability condition for infinitely small $\eps$ should be that of
   the diffusion scheme, that is to say $\Dt\leq \Dx^2/(2\kappa) =
   3\Dx^2\sigma/2$, since $\kappa=1/(3\sigma)$ (see the definition
   after~(\ref{eq-diff})). At the contrary, if the scheme is uniformly
   stable for small $\sigma$ and $\alpha$, the CFL stability condition
   for infinitely small $\sigma$ and $\alpha$ should be that of the
   transport scheme, that is to say $\Dt \leq \eps \Dx$.

The experience with other AP schemes (see the references given in
introduction) suggests that a general condition which is sufficient for
every regime is a combination of the two previous inequalities,
like $\Dt \leq A\max(3\Dx^2\sigma/2, \eps\Dx)$, or $\Dt \leq A(
3\Dx^2\sigma/2 + \eps\Dx)$, and hence is uniform for small $\eps$ and
$\sigma$. This is for instance analytically proved
in~\cite{KU2002} and~\cite{LM2010} for two different schemes. However,
for the UGKS, we were not able to derive such condition so far.

Therefore, in the numerical tests of this paper, we take the empirical condition
$\Dt \leq 0.9(
3\Dx^2\sigma/2 + \eps\Dx)$ that works well for all the regimes we have tested.

Note that if this condition is written with the non-rescaled
(dimensional) variables, we get in fact $\Dt \leq
0.9(3\Dx^2\sigma/(2c)+\Dx/c)$. One could find such a condition still
too restrictive in the transport regime, since $c$ is very large. This is
why implicit schemes are generally preferred in radiation hydrodynamics
for instance. However, a semi-explicit scheme like the UGKS has
the advantage to have a very local stencil, and hence should have a much
better parallel scalability than fully implicit schemes. Indeed, for
large-scale problems run on massively parallel architectures, this
property might make this scheme competitive with fully implicit
schemes. Moreover, there are other applications in which the time step
is limited by some properties of the material, up to a value which is
much lower than the one given by the CFL condition of the UGKS. Finally,
in other problems, like in relativistic hydrodynamics, the time scale
itself is very small, which makes the CFL condition not so
restrictive. See~\cite{MELD_2008} for details.}

%--------------------------------------------------------------------------
\section{Extensions}
\label{sec:extensions}
%--------------------------------------------------------------------------

In this section, we propose some extensions of the previous scheme.

%=========
\subsection{Implicit diffusion }
\label{subsec:ext_imp}
%=========
Our scheme gives for small $\eps$ an explicit scheme for the limit
diffusion equation. This means that at this limit, our scheme requires
the following CFL like condition to be stable: $\Dt\leq
\frac{\Dx^2}{2\mu}$, where $\mu=1/\min(3\sigma)$ is the largest
diffusion coefficient in the domain. When $\Dx$ is small, this
condition is very restrictive. Indeed, diffusion equations are
generally solved by implicit schemes that are free of such a CFL
condition. It is therefore interesting to modify our scheme so as to
recover an implicit scheme in the diffusion limit. Such a modification
is not always trivial (it is not known for schemes
of~\cite{JPT1998,LM2007,klar_sinum_1998,CGLV2008,LS2012}),
 but has already been proposed for some others
(see~\cite{KS2001,Lemou_relaxedAP_2010}). 

For our scheme, we first note that the explicit diffusion term
in~(\ref{eq-num-diffusion-limit}) comes from only one term in the
microscopic numerical flux, namely the fourth term
in~(\ref{eq-flux}): $D_{i+\demi}v^2( \deltax\rhonLipdemi\One_{v>0} +
\deltax\rhonRipdemi\One_{v<0})$. Indeed, the left and right slopes
$\deltax\rhonLipdemi$ and $\deltax\rhonRipdemi$ (defined
by~(\ref{eq-slopes})) depend on the interface value $\rhonipdemi$ and
on the left and right values $\rhoni$ and $\rhonipun$ at time $t_n$. After integration with
respect to $v$, the interface value $\rhonipdemi$ vanishes and we only
get the difference of the $\rhoni$ and $\rhonipun$ {\it at time
  $t_n$}. Consequently, to obtain an implicit diffusion term, that is
to say a difference of left and right values of $\rho$ at time
$t_{n+1}$, it is sufficient to modify the definition of the slopes as
follows: we now set
\begin{equation}  \label{eq-slopes_npun}
\deltax\rhonLipdemi=\frac{\rhonipdemi-\rhonpuni}{\Dx/2} \quad \text{ and
} \quad 
\deltax\rhonRipdemi=\frac{\rhonpunipun-\rhonipdemi}{\Dx/2}.
\end{equation}
Note that the interface value $\rhonipdemi$ is still at $t_n$, so that
the flux is still explicit with respect to $f$. All the other
quantities that are defined our scheme are unchanged.
{\modif  This means that~(\ref{eq-UGKS_rho}) now is an implicit relation
  (that leads to a tridiagonal linear system), while~(\ref{eq-UGKS_f})
  is still explicit with respect to $f$ (in the transport part).}

Now, following the same analysis as in section~\ref{subsec:diffusion_regime},
our scheme converges to the following implicit discrete diffusion
equation:
\begin{equation*}
  \frac{\rhonpuni-\rhoni}{\Dt}   -  \frac{1}{\Dx}
\biggl(
 \frac{1}{3\sigmapdemi}\frac{\rhonpunipun-\rhonpuni}{\Dx}
- 
\frac{1}{3\sigmamdemi}\frac{\rhonpuni-\rhonpunimun}{\Dx}
 \biggr) =- \alpha_i \rhonpuni + G.
\end{equation*}

%=========
\subsection{ More accurate boundary conditions for the diffusion limits}
\label{subsec:ext_CL}
%=========
Here, we propose a modification of the UGKS to obtain the correct
Dirichlet boundary value~(\ref{eq-CLchandra}) in case of a non isotropic
boundary condition. Since
the incorrect value~(\ref{eq-rhondemi_modified}) is required by the
first term $\frac{1}{\eps}\cint{vf_L \One_{v>0}}$ of $\Phi_{\demi}$
in~(\ref{eq-macrosfluxleft}), we propose to replace this term by the
value that gives the correct limit. Indeed, we now impose the relation
\begin{equation}\label{eq-macrosfluxleft_modified}
   \Phi_{\demi}=-\frac{2\cint{v\One_{v<0}}}{\eps}\cint{Wf_L\One_{v>0}}+ \cdots,
% + A_{\demi}\cint{v f^n_1 \One_{v<0}}
% - \frac{1}{4}C_{\demi}\rho^n_{\demi} 
% + \frac{1}{6}D_{\demi} \deltax\rho^{nR}_{\demi}\One_{v<0}
% - \frac{1}{4} E_{\demi}G\One_{v<0}.
\end{equation}
where $\cdots$ stands for the other terms of~(\ref{eq-macrosfluxleft})
that are unchanged.  Moreover, we set
\begin{equation}\label{eq-BC_modified}
\rho^n_{\demi}= 2\cint{W f_L \One_{v>0}}
\end{equation}
in~(\ref{eq-macrosfluxleft_modified}) and~(\ref{eq-fluxleft}). These
definitions ensure that $\Phi_{\demi}$ is uniformly bounded with
respect to $\eps$, and moreover, the discrete conservation
law~(\ref{eq-conslaws}) in the first cell $i=1$ tends
to~(\ref{eq-diff1}) with the correct Dirichlet boundary data $2\cint{W f_L \One_{v>0}}$.

However, contrary to the definition of $\rho^n_{\demi}$ suggested in
section~\ref{subsec:BC}, this definition has an influence in the free
transport regime. Indeed, while $\phi_{\demi}$ tends to the correct first
order upwind boundary flux, the macroscopic boundary flux
$\Phi_{\demi}$ tends to $\frac{1}{\eps}\cint{Wf_L\One_{v>0}+ v f^n_1\One_{v<0}}$,
which is not the correct flux for the free transport regime (it should
be $\frac{1}{\eps}\cint{vf_L\One_{v>0}+ v f^n_1\One_{v<0}}$).

Finally, we can obtain the two correct limits by using a blended
modification of the first term of $\Phi_{\demi}$
in~(\ref{eq-macrosfluxleft}). Indeed, we set
\begin{equation}\label{eq-macrosfluxleft_remodified}
   \Phi_{\demi}=\frac{1}{\eps}\Bigl\langle \bigl[(1-\theta(\nu_{\demi}))v 
+ \theta(\nu_{\demi})\times(-2\cint{v\One_{v<0}})W\bigr]f_L\One_{v>0}\Bigr\rangle 
+ \cdots,
% + A_{\demi}\cint{v f^n_1 \One_{v<0}}
% - \frac{1}{4}C_{\demi}\rho^n_{\demi} 
% + \frac{1}{6}D_{\demi} \deltax\rho^{nR}_{\demi}\One_{v<0}
% - \frac{1}{4} E_{\demi}G\One_{v<0}.
\end{equation}
where again $\cdots$ stands for the other terms
of~(\ref{eq-macrosfluxleft}). We set 
\begin{equation} \label{eq-BC_blended}
\rho^n_{\demi}= 2\Bigl\langle \bigl[((1-\theta(\nu_{\demi}))v 
+ \theta(\nu_{\demi})W\bigr] f_L \One_{v>0}\Bigr\rangle 
\end{equation}
in~(\ref{eq-macrosfluxleft_modified}) and~(\ref{eq-fluxleft}), and the
blending parameter is $\theta(\nu)=1-\exp(-\nu\Dt)$. This parameter is
such that $\theta(\nu)$ tends to 1 for small $\eps$ and tends to 0 for
small $\sigma$ and $\alpha$.  This definition implies that in the
diffusion limit, we get indeed the correct Dirichlet boundary
condition $2\cint{W f_L \One_{v>0}}$, and that in the free transport
regime, the microscopic and macroscopic boundary fluxes tend to the
correct values. Note that these properties are illustrated in
section~\ref{sec:res}.

%=========
\subsection{ General linear collision operators}
%=========
In this section, we propose a way to adapt the UGKS approach to a
general linear Boltzmann operator. Namely, we consider the equation
\begin{equation}  \label{eq-linearbtz}
\eps \dt f + v \dx f = \frac{1}{\eps}Lf ,
\end{equation}
where the collision operator now is
\begin{equation}  \label{eq-L}
Lf=\int_{-1}^1 k(v,v')(f(v')-f(v))\, dv'.
\end{equation}
We ignore the absorption and source terms here, since this does not
change our analysis. For such a model, it is well known that the density
of $f$ satisfies in the small $\eps$ limit the following diffusion
equation
\begin{equation}  \label{eq-diffL}
\dt \rho +\dx\kappa\dx \rho=0,
\end{equation}
with $\kappa=\cint{vL^{-1}v}$, and where $L^{-1}$ is the
pseudo-inverse of $L$. This operator is defined for functions with
zero average by the following: for any $\phi$ such that
$\cint{\phi}=0$, $\psi=L^{-1}\phi$ is the unique solution of
$L\psi=\phi$ such that $\cint{\psi}=0$.

Since in the UGKS approach the relaxation form of the collision
operator is strongly used, a natural idea is to write $Lf$ in the
gain-loss form: $Lf=L_+f-\frac{1}{\tau(v)}f$, where $1/\tau(v)=\int_{-1}^1 k(v,v')\,
dv'$. Then we can try to apply the previous strategy in which the gain
term $L_+f$ plays the role of $\rho$. However, we observed that this
strategy fails completely, since the resulting scheme cannot capture the correct
diffusion limit. 
% More precisely, this scheme is
% \begin{equation}\label{eq-fvL}
% \frac{\fnpuni-\fni}{\Dt}  + \frac{1}{\Dx}\left(
%   \phipdemi - \phimdemi
%   \right) 
%  = \frac{1}{\eps^2}L\fnpuni,
% \end{equation}
% where the numerical fluxes are defined by
% $\phipdemi=\frac{1}{\eps\Dt}\int_{t_n}^{\tnpun} v f(t,\xipdemi,v)  \,
% dt $, in which $f(t,\xipdemi,v)$ is approximated by the
% exponential-integral representation of the solution of~(\ref{eq-linearbtz}):
% \begin{equation*}
% \begin{split}
%   f(t,\xipdemi,v) & \approx   e^{-\nu(t-\tn)}f(\tn,\xipdemi-\frac{v}{\eps}(t-t_n)) \\
% & \quad + \frac{1}{\eps^2}\int_{t_n}^{t} e^{-\nu(t-\tn)}L_+f(s,\xipdemi-\frac{v}{\eps}(t-s))  \, ds,
% \end{split}
% \end{equation*}
% where $\nu(v)=1/\eps^2\tau(v)$. After some algebra, the fluxes can
% be computed analytically, and an asymptotic analysis gives the
% following result: in the small $\eps$ limit, $\rho^n$ satisfies a
% relation which is not consistent with the diffusion limit. In order to
% shorten the article, we do not
% give the details of this assertion here. Note that another problem is that
% scheme~(\ref{eq-fvL}) requires the inversion of $I-\Dt/\eps^2L$, since
% using the conservation law is not sufficient here, which can be
% expensive.

Then we propose a different and simple way to capture the
correct asymptotic limit by using the penalization technique of
Filbet and Jin~\cite{FJ2010}. First, the collision operator is written
as
\begin{equation}  \label{eq-penL}
Lf=(Lf-\theta Rf) +\theta Rf,
\end{equation}
where $Rf=\rho-f$ (with $\rho=\cint{f}$) is the corresponding
isotropic operator, and where $\theta$ is a parameter adjusted to capture the
correct diffusion coefficient. Therefore,
equation~(\ref{eq-linearbtz}) can be rewritten as
\begin{equation}  \label{eq-penbtz}
\eps \dt f + v \dx f = \frac{\theta}{\eps}(\rho-f) + \eps\tilde{G},
\end{equation}
where $\tilde{G}=(Lf-\theta Rf)/\eps^2$ is considered as a time and
velocity dependent source. Then our previous derivation can be readily
applied to this equation: we replace $G$ by $\tilde{G}$, $\sigma$
by $\theta$, and $\alpha$ by 0 in the different steps of
sections~\ref{subsec:FV} to~\ref{subsec:summary} to get the following
scheme:
\begin{align*}
&  \frac{\rhonpuni-\rhoni}{\Dt}  + \frac{1}{\Dx}\left(
  \Phipdemi - \Phimdemi
  \right) 
 = 0, \\
& \frac{\fnpuni-\fni}{\Dt}  + \frac{1}{\Dx}\left(
  \phipdemi - \phimdemi
  \right) 
 = \frac{\theta}{\eps^2}(\rhonpuni-\fnpuni) +
 \frac{1}{\eps^2}(Lf^n_i-\theta R\fni).
\end{align*}
The numerical fluxes are given by
relations~(\ref{eq-flux}--\ref{eq-coeff_flux}), in which we just have
to replace $G$ by $\tilde{G}$, $\sigma$
by $\theta$, and $\alpha$ by 0 (note that $\nu$ must be modified
accordingly in these relations).

Consequently, our analysis detailed in section~\ref{sec:analysis}
directly applies to this scheme, and we obtain that it gives
for small $\eps$ this limit diffusion scheme
\begin{equation*}
  \frac{\rhonpuni-\rhoni}{\Dt}   -  \frac{1}{\Dx}
\biggl(
 -\frac{\cint{v^2}}{\theta}\frac{\rhonipun-\rhoni}{\Dx}
+ 
\frac{\cint{v^2}}{\theta}\frac{\rhoni-\rhonimun}{\Dx}
 \biggr) =0,
\end{equation*}
which is a consistent approximation of the diffusion
limit~(\ref{eq-diffL}), if and only if the parameter $\theta$ is set
to
\begin{equation}\label{eq-deftheta} 
  \theta=-\frac{\cint{v^2}}{\cint{vL^{-1}v}}.
\end{equation} 
This proves that this "penalized" UGKS is AP for the limit $\eps
\rightarrow 0$.

Note that to use this scheme, we just have to compute $L^{-1}v$, which
has to be done only once. This computation can be made analytically
(for simple scattering kernels) or numerically. Also note that an
approach with similar aspects has been proposed
in~\cite{Lemou_relaxedAP_2010}.

However, this AP property is true only if the scheme is uniformly
stable. Since this property is difficult to obtain for the complete
scheme, we restrict to the space homogeneous problem, and we show that
this leads to a non trivial restriction on the collision operator.
\begin{proposition}
  The penalized scheme
\begin{equation*}
   \frac{f^{n+1}-f^n}{\Dt} 
 = \frac{\theta}{\eps^2}(\rhonpun-f^{n+1}) +
 \frac{1}{\eps^2}(Lf^n-\theta R f^n)
 \end{equation*}
 that approximates the homogeneous equation $\dt f = Lf$ is absolutely
 stable if $\Dt(k_M-\theta)\leq \eps^2$, where $k_M=\max_{v,v'} k$. This
 stability is uniform with respect to $\eps$ if $\theta\geq k_M$.
\end{proposition}
This proposition can be easily proved by using the fact that if the scattering
kernel is bounded ($0<k_m\leq k(v,v')\leq k_M$), then $L$ is a non
positive self-adjoint operator, and that its eigenvalues are all
bounded in absolute value by $2k_M$.

 Note that with $\theta$ defined by~(\ref{eq-deftheta}), this
 restriction reads $-\cint{vL^{-1}v}\leq \cint{v^2}/k_M$. If this
 inequality is not satisfied, then $\Dt$ must be smaller and smaller as $\eps$
 decreases, and hence the scheme cannot be AP. This property will be
 studied in detail in a forthcoming work for realistic functions $k$
 (like the Henyey-Greenstein function).

%--------------------------------------------------------------------------
\section{Numerical results}
\label{sec:res}
%--------------------------------------------------------------------------

Almost all the test cases presented here are taken from
references~\cite{klar_sinum_1998,JPT1998}. Comparisons of several
existing AP schemes with these test cases can be found
in~\cite{LM2007}. Depending on the regime, we compare UGKS to
a standard upwind explicit discretization of~(\ref{eq-onegroup}) or to the
explicit discretization of the diffusion limit~(\ref{eq-diff}). These
reference solutions are obtained after a mesh convergence
study: the number of points is sufficiently large to consider that the
scheme has converged to the exact solution. Generally, the time step
for the UGKS is taken as $\Dt=cfl\max(\eps\Dx,3\Dx^2\sigma/2)$, where $cfl=0.9$.

%=========
   \subsection{Various test cases with isotropic boundary conditions}
   \label{subsec:isotropic}
%=========

\begin{figure}
\centerline{\includegraphics[width=0.6\textwidth]{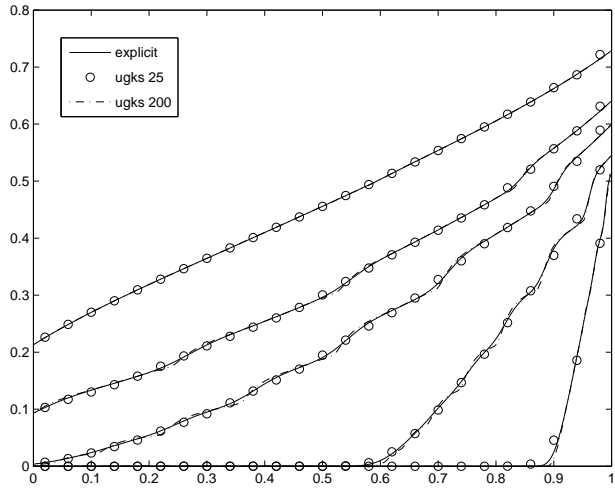}}
\caption{\label{fig:KLAR_E1_ugks} Kinetic regime (example 1): comparison between
  a reference solution and the UGKS (25 and 200 grid points). Results at times $t=0.1$, $0.4$,
  $1.0$, $1.6$ and $4$ ($\eps=1$).}
\end{figure}

\paragraph{Example 1}
Kinetic regime:
\begin{align*}
& x\in[0,1] , \qquad  f_L(v)=0 , \qquad  f_R(v)=1, \\
& \sigma=1 , \qquad  \alpha=0 , \qquad  G=0 , \qquad  \eps=1.  
\end{align*}
The results are plotted at times $t=0.1$, $0.4$, $1.0$, $1.6$, and
$4$. We use 25 and 200 points for UGKS. The reference solution is
obtained with 1000 points. In figure~\ref{fig:KLAR_E1_ugks}, we observe that
the UGKS is very close to the reference solution, even with the
coarse discretization for
short times $t=0.1$ and $t=0.4$ (which is better than the AP schemes
compared in~\cite{LM2007}).

\begin{figure}
\centerline{\includegraphics[width=0.6\textwidth]{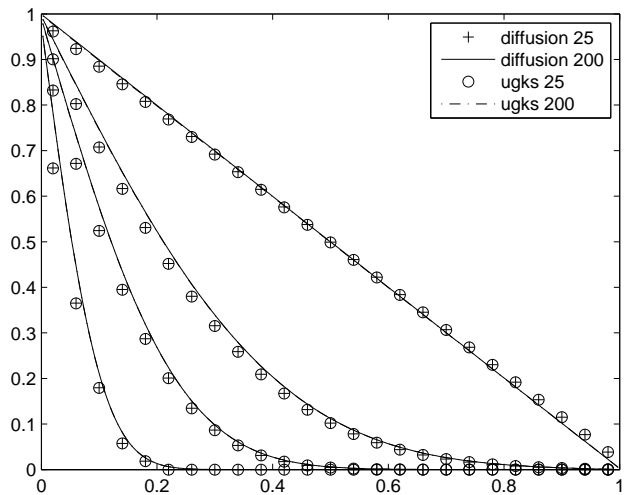}}
\caption{\label{fig:JPT1_ugks} Diffusion regime (example 2): comparison between diffusion
  solution and the UGKS (25 and 200 grid points). Results at
  times $t=0.01$, $0.05$, $0.15$ and $2$ ($\eps=10^{-8}$).}
\end{figure}

\paragraph{Example 2}
Diffusion regime:
\begin{align*}
& x\in[0,1] , \qquad  f_L(v)=1 , \qquad  f_R(v)=0, \\
& \sigma=1 , \qquad  \alpha=0 , \qquad  G=0 , \qquad  \eps=10^{-8}.  
\end{align*}
The results are plotted at times $t=0.01$, $0.05$, $0.15$ and $2$. We
use 25 and 200 points for the UGKS. Here, the reference solution is obtained
with the explicit discretization of the diffusion equation, since the
kinetic equation cannot be solved by a the standard upwind scheme with such a small $\eps$.  In figure~\ref{fig:JPT1_ugks}, we
see that the UGKS and the diffusion solution are almost indistinguishable at any times for both coarse and fine discretizations.

\begin{figure}[p]
\centerline{\includegraphics[width=0.6\textwidth]{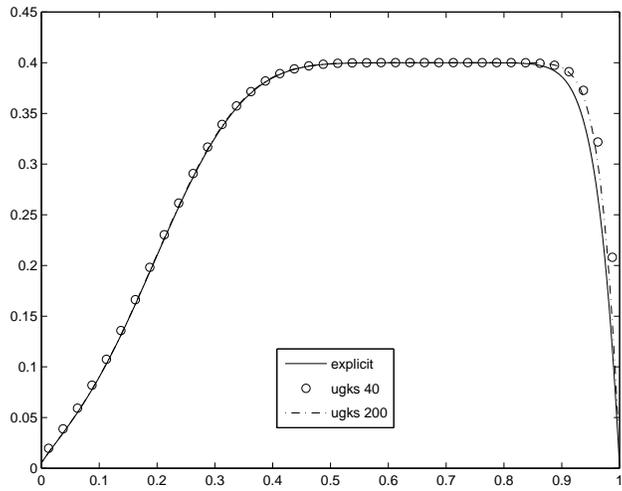}}
\caption{\label{fig:JPT2_ugks}Intermediate regime with a variable
  scattering frequency and a source term (example 3): comparison between a
  reference solution and the UGKS (40 and 200 grid points).}
\end{figure}

\paragraph{Example 3} 
Intermediate regime with a variable scattering frequency and a source term:
\begin{align*}
& x\in[0,1] , \qquad  f_L(v)=0 , \qquad  f_R(v)=0, \\
& \sigma=1+(10x)^2 , \qquad  \alpha=0 , \qquad  G=1 , \qquad  \eps=10^{-2}.  
\end{align*}
The results are plotted at times $t=0.4$ with 40 and 200 points in
figure~\ref{fig:JPT2_ugks}. The
reference solution is obtained with the explicit
scheme using 20 000 points. We observe that the UGKS provides results
that are very close to the reference solution, like the AP schemes
compared in~\cite{LM2007}.

\begin{figure}[p]
\centerline{\includegraphics[width=0.6\textwidth]{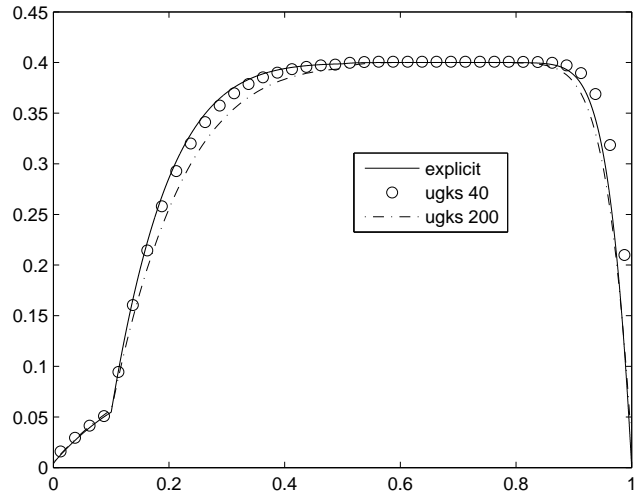}}
\centerline{\includegraphics[width=0.6\textwidth]{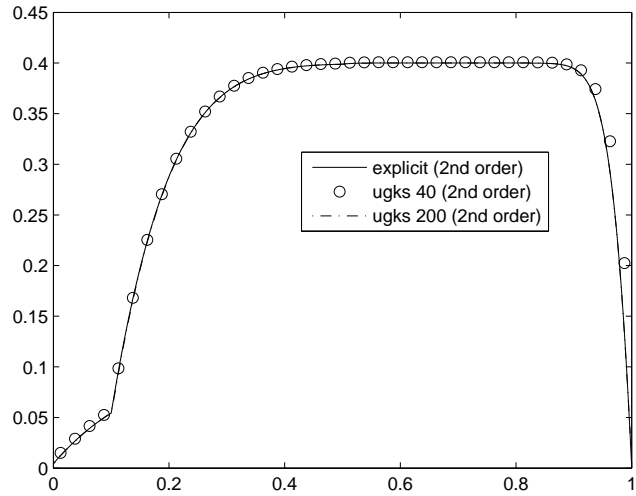}}
\caption{\label{fig:cdiscont}{\modif Intermediate regime with a discontinuous
  scattering frequency and a source term (example 4). Comparison between a
  reference solution and the UGKS (40 and 200 grid points): (top)
  first order schemes, (bottom) second order schemes.}
}
\end{figure}

{\modif 
\paragraph{Example 4}
Intermediate regime with a discontinuous scattering frequency and a
source term. We take the following discontinuous values of $\sigma$: 
\begin{align*}
 x & \in[0,0.1] , &  x &\in [0.1,0.5]  ,&   x& \in[0.5,1], \\
 \sigma & =1 , &  \sigma & =10 , &  \sigma& =100,
\end{align*}
while all the other parameters are like in example 3. In
figure~\ref{fig:cdiscont} (top), we show the behavior of the UGKS for
coarse and thin meshes, as compared to a reference solution obtained
by a fully explicit scheme with 20\,000 cells.  We observe that the
UGKS gives results that are close to the reference solution. However,
the result with the thin mesh is not as good as expected in the small
to moderate scattering region ($x<0.4$). This is in fact due to the
numerical dissipation induced by the first order reconstruction of
$f(t_n)$ around the cell interface (see~(\ref{eq-piecewisef})), since
for small scattering, this gives a standard first order upwind scheme,
which is known to be very diffusive.

This is confirmed by the following modification. We use a second order
reconstruction with slope limiters to replace~(\ref{eq-piecewisef}) by
\begin{equation*}
f(t_n,x,v)=\left\lbrace
\begin{split}
& \fni+\frac{\Delta x}{2}\delta f^n_i \quad \text{ if } \quad x<\xipdemi \\
& \fnipun-\frac{\Delta x}{2}\delta f^n_{i+1} \quad \text{ if } \quad x>\xipdemi,
\end{split}
\right.
\end{equation*}
where $\delta f^n_i$ is the slope computed with the MC limiter:
$\delta f^n_i=\text{minmod}(\frac{f^{n}_{i+1}-f^{n}_{i-1}}{2\Delta
  x},\theta\frac{f^{n}_{i}-f^{n}_{i}}{\Delta x},$ \\
$\theta\frac{f^{n}_{i+1}-f^{n}_{i}}{\Delta x})$ and
$\theta=1.5$. Then, following the construction explained in
section~\ref{sec:ugks}, we again obtain
scheme~(\ref{eq-UGKS_rho})--(\ref{eq-UGKS_f}), in which the additional
term $B_{i+\demi}v (v\delta f^n_{i}\One{v>0} + v\delta
f^n_{i+1}\One{v<0})$ must be added to the numerical flux given
in~(10), with
\begin{equation*}
 B(\Dt, \eps, \sigma,\alpha) =\frac{1}{\Dt \eps^2 \nu}
                                \Bigl(\Dt e^{-\nu \Dt} -
  \frac{1}{\nu}\bigl(1-e^{-\nu \Dt}\bigr)\Bigr).
\end{equation*}
We also recompute the reference solution, still with 20\,000 points,
but with the same second order upwind reconstruction. We observe in
figure~\ref{fig:cdiscont} (bottom) that the UGKS results are now much
closer to the reference solution, and that there is no visible
difference between this solution and the UGKS with the thin mesh.
}

%=========
   \subsection{Comparison of different numerical boundary conditions}
   \label{subsec:res_bc}
%=========

In this test, we compare the different numerical boundary
conditions (BC) proposed for the UGKS in section~\ref{subsec:BC}
and~\ref{subsec:ext_CL}. The BC given
in~(\ref{eq-rhondemi_modified}) that ensures the stability of the
scheme is called {\it stabilized BC}. We call {\it corrected BC} the
one described
in~(\ref{eq-macrosfluxleft_modified},\ref{eq-BC_modified}) that
gives a correct boundary value in the diffusion limit. Finally,
the {\it blended BC} is that defined 
in~(\ref{eq-macrosfluxleft_remodified},\ref{eq-BC_blended}).

\paragraph{Example 5}
First, we consider an intermediate regime with a non-isotropic boundary
condition that generates a boundary layer of size $L=0.01$ at the left
boundary:
\begin{align*}
& x\in[0,1] , \qquad  f_L(v)=v , \qquad  f_R(v)=0, \\
& \sigma=1, \qquad  \alpha=0 , \qquad  G=0 , \qquad  \eps=10^{-2}
\end{align*}

\begin{figure}
\centerline{\includegraphics[width=0.6\textwidth]{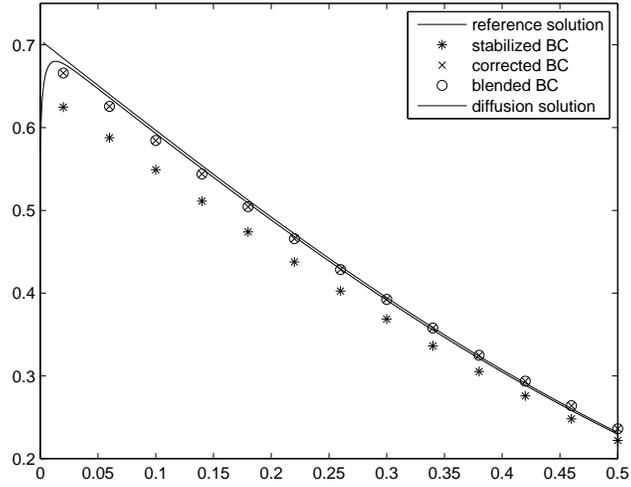}}
\centerline{\includegraphics[width=0.6\textwidth]{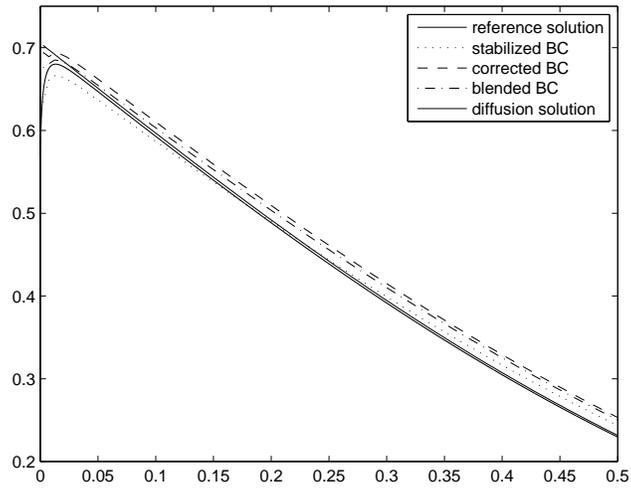}}
\caption{\label{fig:klar_E4_eps2_sigma1_alpha0}Boundary layer problem
  (example 5):
  stabilized, corrected, and blended boundary conditions for the UGKS
  are compared to the explicit scheme (solid line with boundary layer)
  and the diffusion solution (solid straight line), $\eps=10^{-2}$. }
\end{figure}

The results are plotted at times t = 0.4 with 25 and 200 points in
figure~\ref{fig:klar_E4_eps2_sigma1_alpha0} (see~\cite{LM2007} for a
comparison of other AP schemes on the same test case). The reference
solution is obtained with the explicit scheme using 20 000 points. We
also consider the diffusion scheme with 25 points, with a left
boundary condition $\rho_L=17/24$ computed by formula~(\ref{eq-CLchandra})
with the approximation $W(v)=3/2v^2+v$. With the coarse
discretization, the boundary layer is of course not resolved. The
stabilized BC is rather different from the reference solution in all
the domain. At the contrary, we observe that the corrected and blended
BC are very close to the reference solution inside the domain (far
from the boundary).

With the fine discretization, the boundary layer is resolved:
surprisingly, the corrected and blended BC are less close to the reference
solution, while the stabilized version is now closer, even in the
boundary layer. This fact is rather different from the other AP
schemes already used for the this test
(see~\cite{LM2007}). We do not understand well the reason so far, but
it might be linked to the fact that the coefficients $A,C,D$ of the
UGKS depend in an intricate way on the numerical parameter $\Dt$, $\Dx$
and $\eps$: for instance for a given small $\eps$, $A$ and $D$ are
closer to their limit value $0$ and $-1/\sigma$ for large $\Dt$ and
$\Dx$ than for smaller values. See a similar observation
in~\cite{TX2006} about the high accuracy of the gas kinetic scheme (a
macroscopic version of the UGKS) for under resolved cases.

\begin{figure}
\centerline{\includegraphics[width=0.6\textwidth]{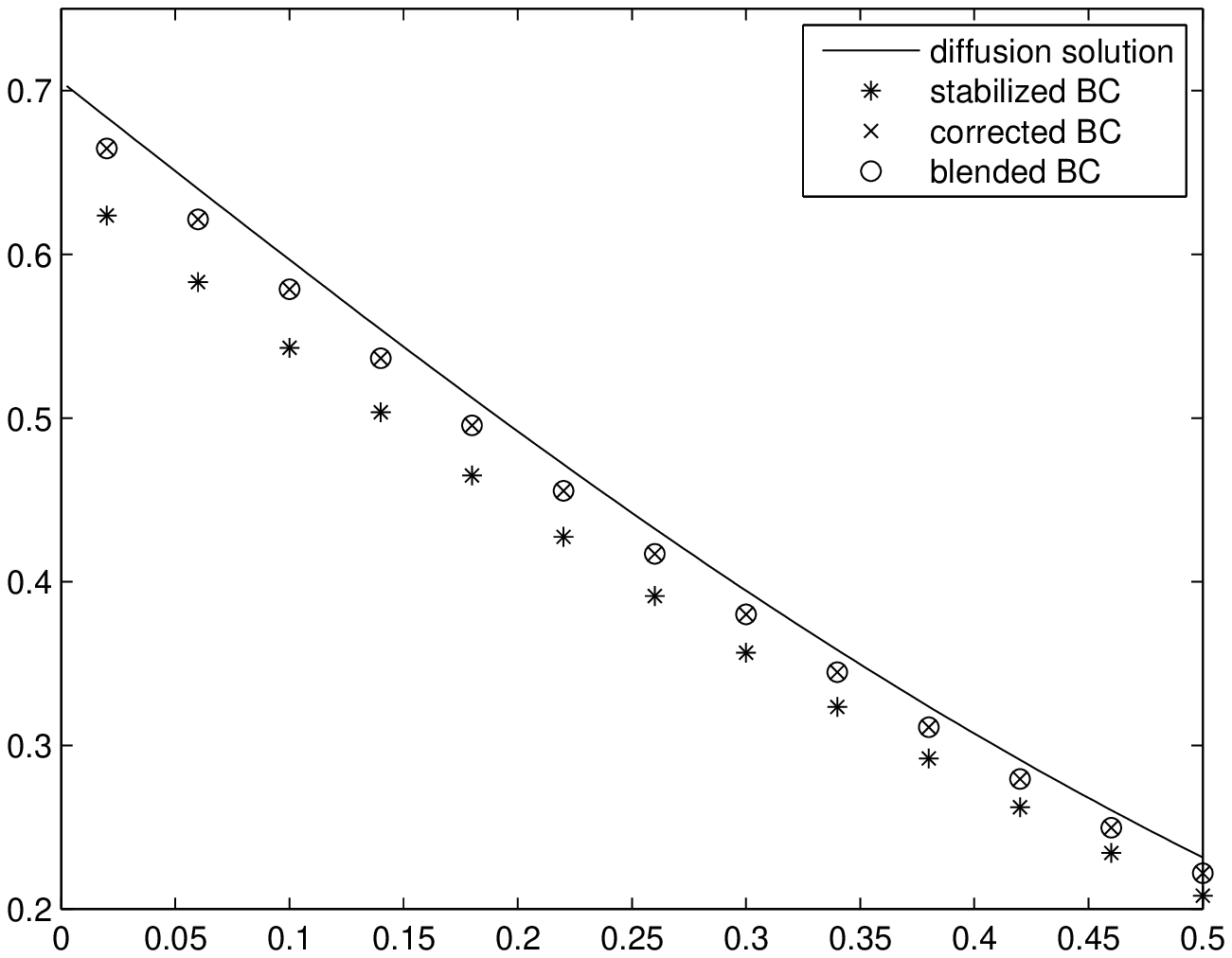}}
\centerline{\includegraphics[width=0.6\textwidth]{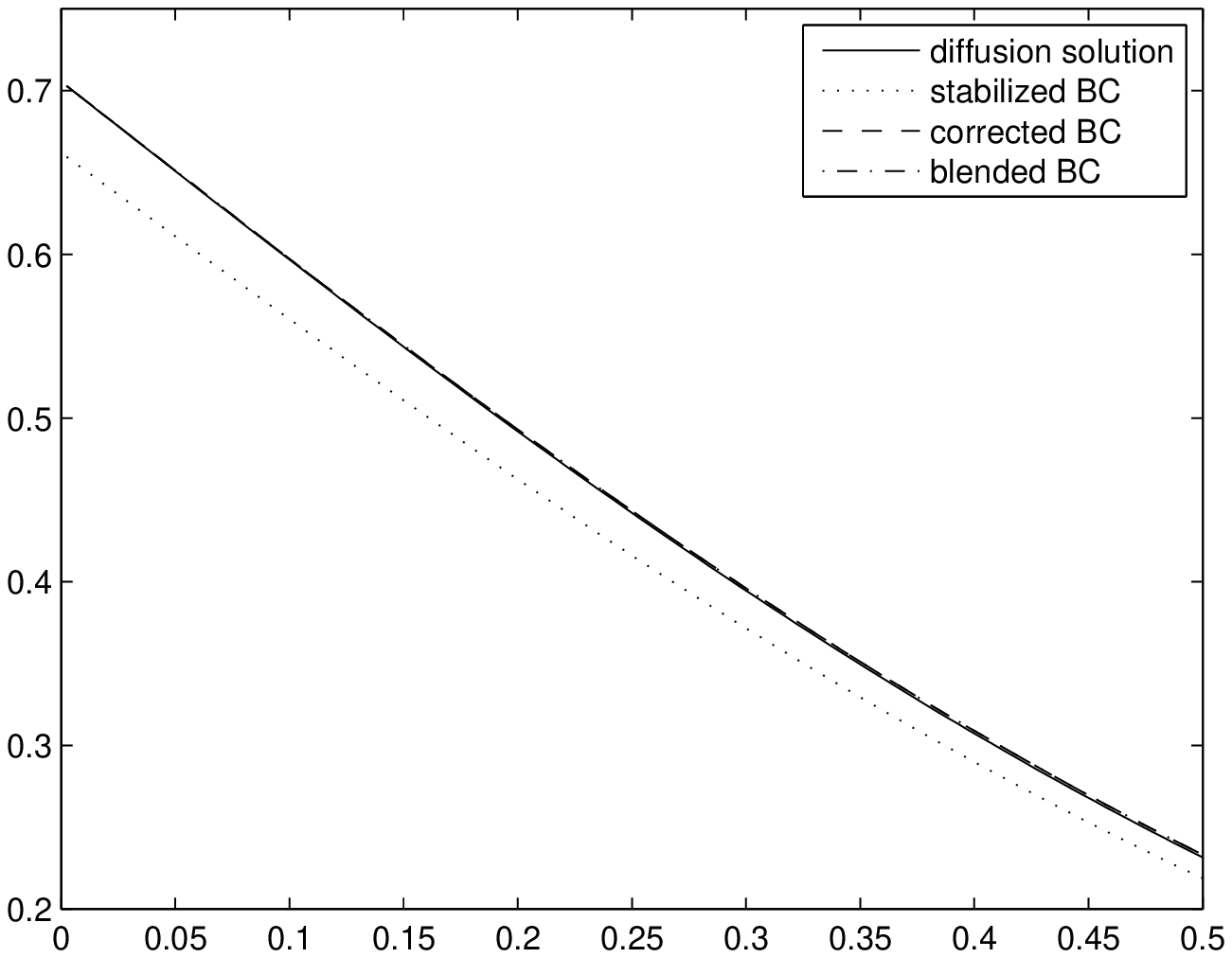}}
\caption{\label{fig:klar_E4_eps4_sigma1_alpha0}Boundary layer problem
  (example 6):
  stabilized, corrected, and blended boundary conditions for the UGKS
  are compared to the diffusion solution, $\eps=10^{-4}$.}
\end{figure}

\paragraph{Example 6}
Then we take $\eps=10^{-4}$ to be in a fully diffusion
regime. The boundary layer is here very small and cannot be
captured. The reference solution then is the solution of the diffusion
equation, computed with $2000$ cells. The results are plotted at
$t=0.4$ in figure~\ref{fig:klar_E4_eps4_sigma1_alpha0}. The comparison is almost the
same as with $\eps=10^{-2}$ for the coarse mesh: the stabilized
BC gives a result which is not very accurate, while, as
expected, the corrected and blended BC (which give
the same results) are much closer to the diffusion solution. However,
with the fine mesh, the corrected and blended BC now give
no difference with the diffusion solution, while the stabilized
BC is not accurate.

\begin{figure}
\centerline{\includegraphics[width=0.6\textwidth]{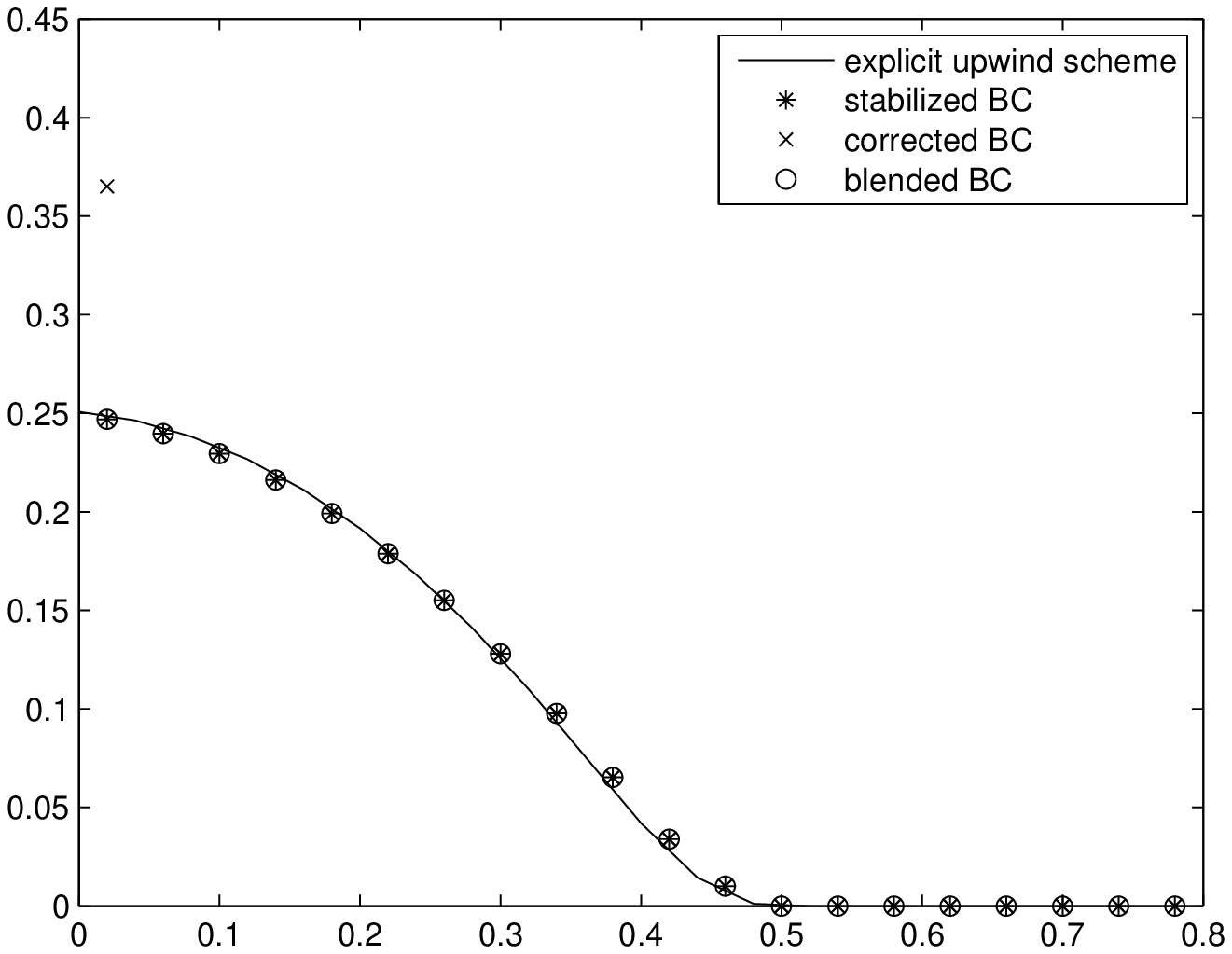}}
\centerline{\includegraphics[width=0.6\textwidth]{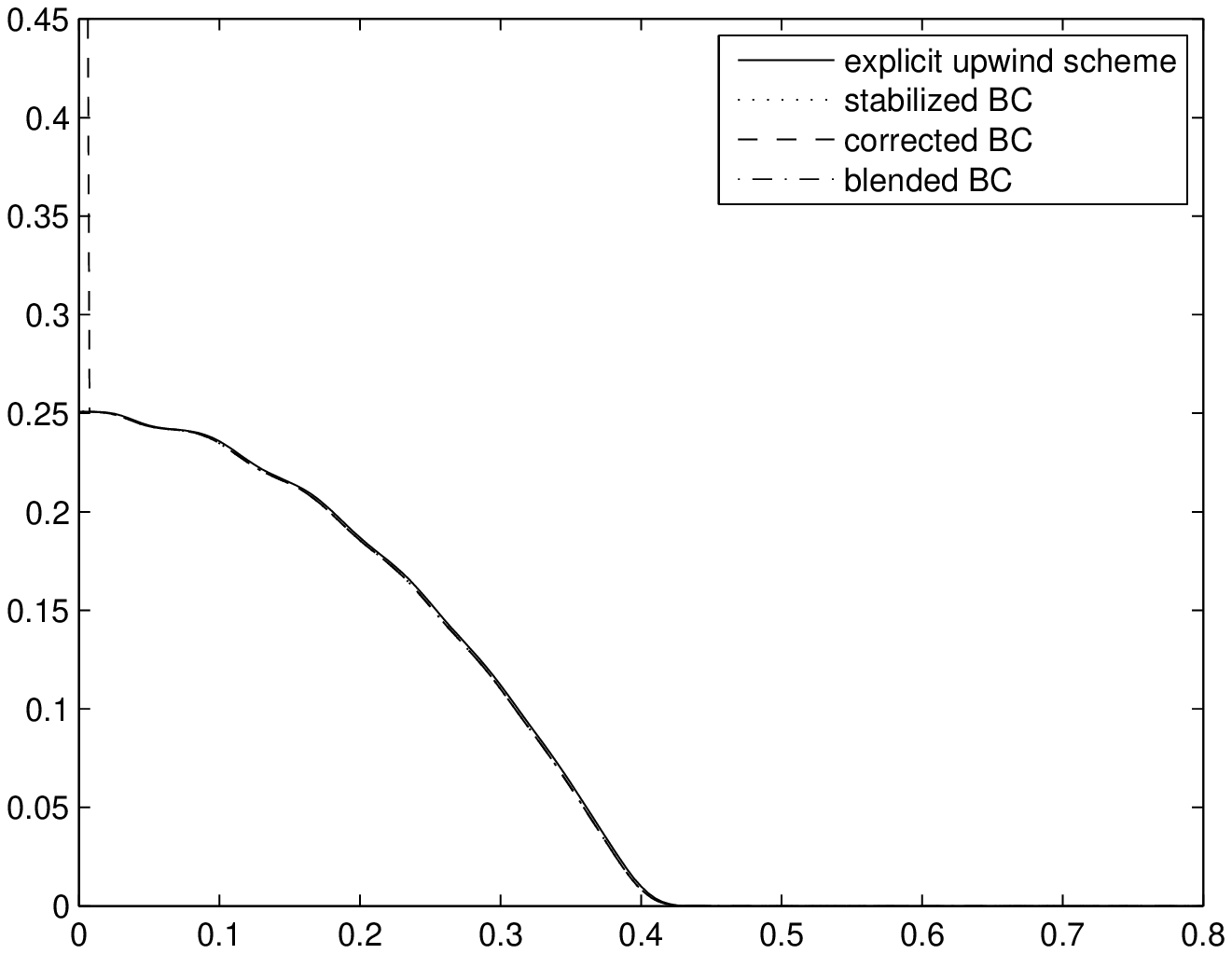}}
\caption{\label{fig:klar_E4_eps0_sigma0_alpha0}Free transport regime
  (example 7): 
  stabilized, corrected, and blended boundary conditions for the UGKS
  are compared to the upwind explicit scheme, $\eps=1$, $\sigma=0$.}
\end{figure}

\paragraph{Example 7}
Finally, in order to test our blended BC in all the
regimes, we take $\eps=1$ and $\sigma=1$, which is a purely transport
regime (no collision). 
The results are plotted at times $t = 0.4$ with 25 and 200 points in
figure~\ref{fig:klar_E4_eps0_sigma0_alpha0}. Here, the comparison with
a reference solution is difficult: indeed, it is well known that
discrete ordinates methods in free transport regimes require many
velocity points to be accurate. With a small number of velocities
(like the 16 Gauss points we use here), the results show the well known ``ray
effect'' that looks like several plateaus in the solution, and this
phenomenon is stronger when the space mesh is refined. Consequently,
we think it more interesting to compare the UGKS to a standard upwind
explicit scheme with the same space resolution. As expected, we observe in
figure~\ref{fig:klar_E4_eps0_sigma0_alpha0} that both stabilized and
blended BCs are very close to the standard scheme (for
both coarse and refine space grids). At the contrary, the corrected
BC (which is not consistent in the free transport
regime) gives results that are far from the other schemes close to the
boundary: it gives a numerical boundary layer at the left boundary
that makes the solution much too large.

%=========
   \subsection{Implicit diffusion}
   \label{subsec:res_imp}
%=========
Here, we illustrate the properties of the UGKS modified to recover an
implicit scheme in the diffusion limit (see
section~\ref{subsec:ext_imp}). For this scheme, we take a time step
defined by $\Dt=\max(0.9\eps\Dx,cfl\Dx)$--where $cfl$ depends
on the problem (between $0.9$ and $0.1$)--such that we get
$\Dt=cfl\Dx$ in the diffusion limit instead of $\Dt=cfl\Dx^2/2\mu$.

\begin{figure}
\centerline{\includegraphics[width=0.6\textwidth]{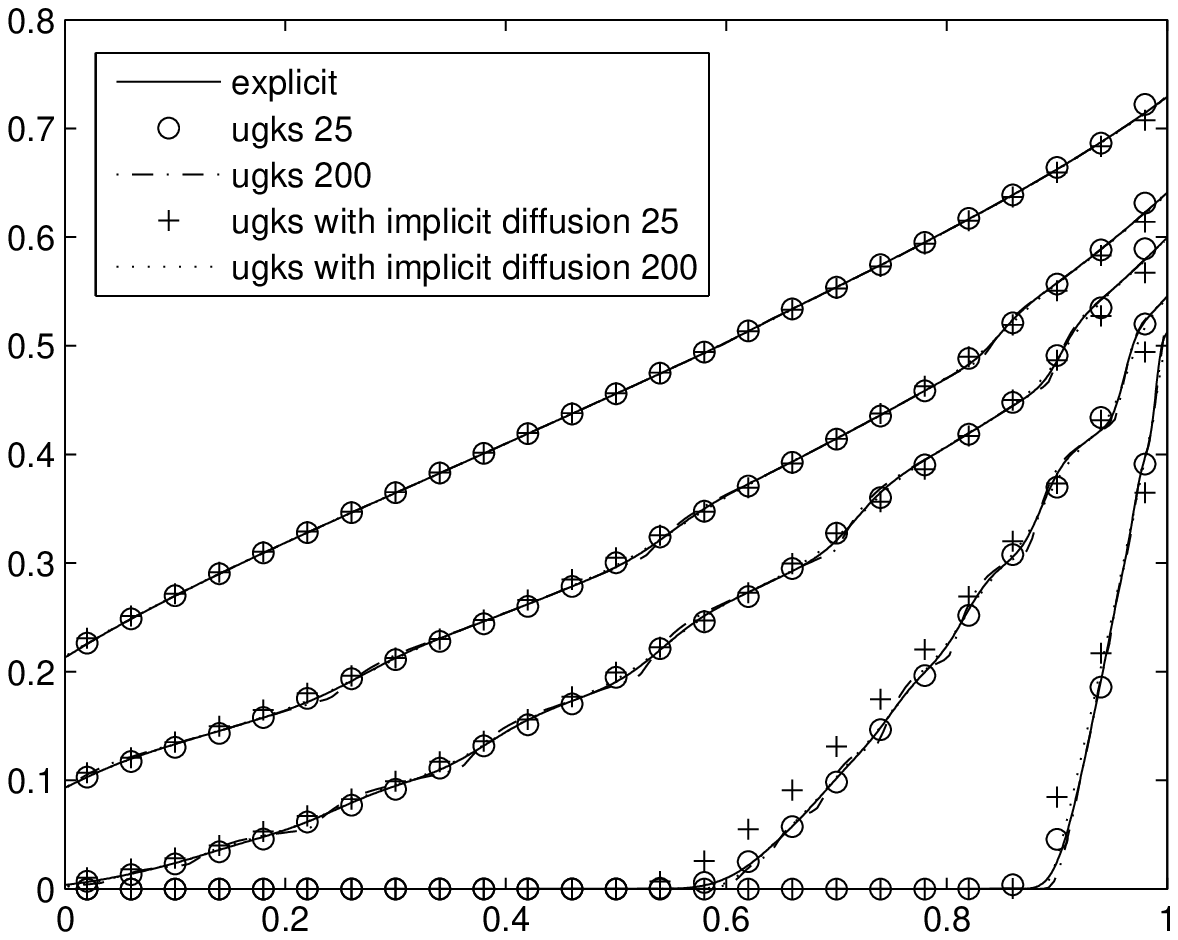}}
\centerline{\includegraphics[width=0.6\textwidth]{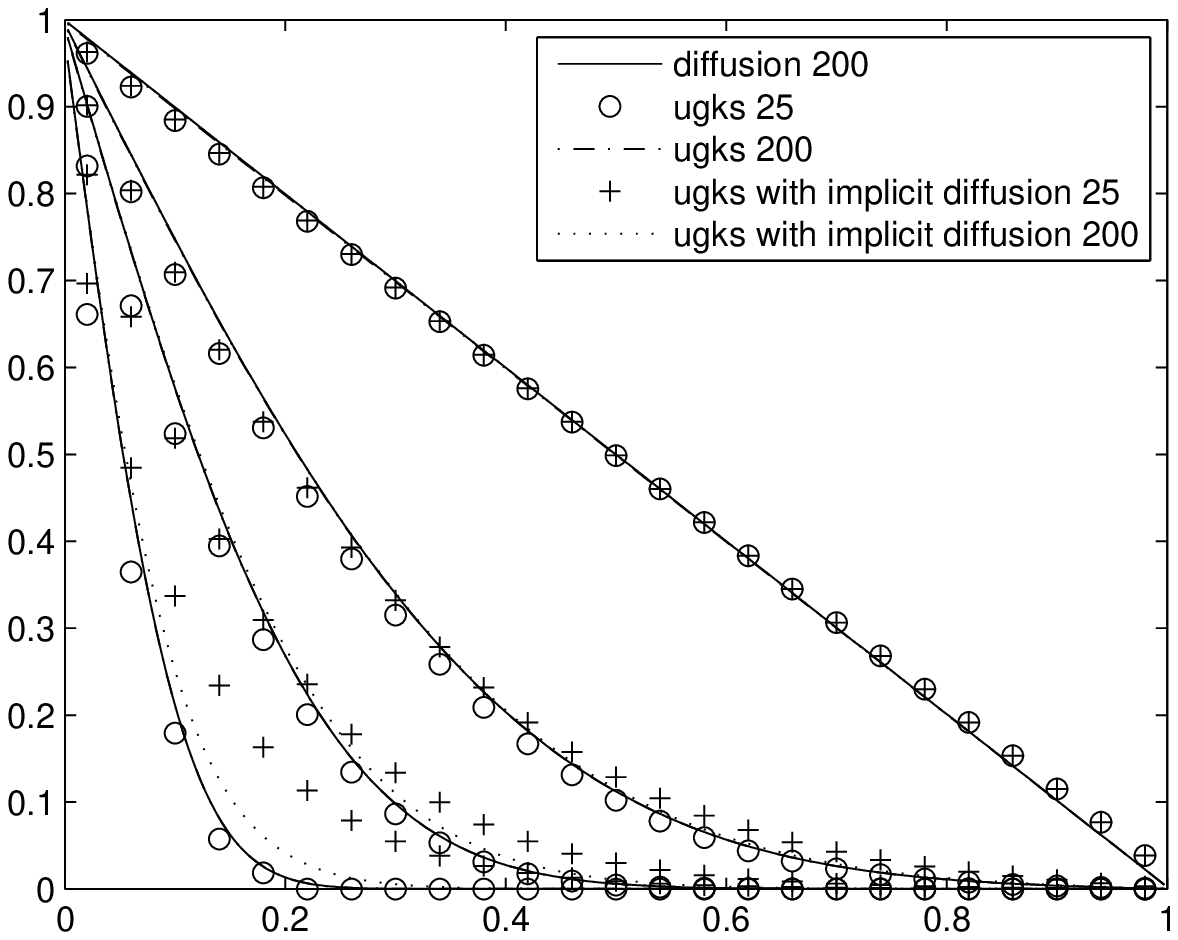}}
\caption{\label{fig:diff_imp_E1_JPT1}Comparison of the UGKS with/without
  implicit diffusion with a reference solution for unsteady cases. Data of of examples
  1 (kinetic regime, top)
  and 2 (diffusion regime, bottom).}
\end{figure}

In figure~\ref{fig:diff_imp_E1_JPT1}, we compare this modified UGKS
(denoted by UGKS-ID in the following) to
the non-modified UGKS (denoted by UGKS-ED) and to the reference solutions for two unsteady
cases. The data are those of example 1 (with a kinetic regime) and 2
(with a diffusion regime) shown in section~\ref{subsec:isotropic}. For
example 1, we are in a kinetic regime where our modification is
useless. However, with the fine mesh, the two UGKS give results that
are very close. With the coarse mesh, there are some differences, but
still quite small. For example 2, we are in a diffusion regime, where
the UGKS-ID allows to take a time step which does not respect
the parabolic CFL $\Dt<\frac{\Dx^2}{2\mu}$. Indeed, in the case of the
fine mesh, while the UGKS-ID
requires a time step of $3.37 \, 10^{-5}$, the UGKS-ED
requires a time step of only $4.5\, 10^{-2}$, which is $133$ times
larger. In addition, we observe that both schemes give results that
are very close.
Of course, the results are different for the small times (the schemes are first
order in time, hence a larger time step induces a larger numerical
error), but we see that for longer times, both schemes give very close
results. As expected, if the time step of the UGKS-ID is
decreased, the difference with the UKGS-ED reduces too.

\begin{figure}
\centerline{\includegraphics[width=0.6\textwidth]{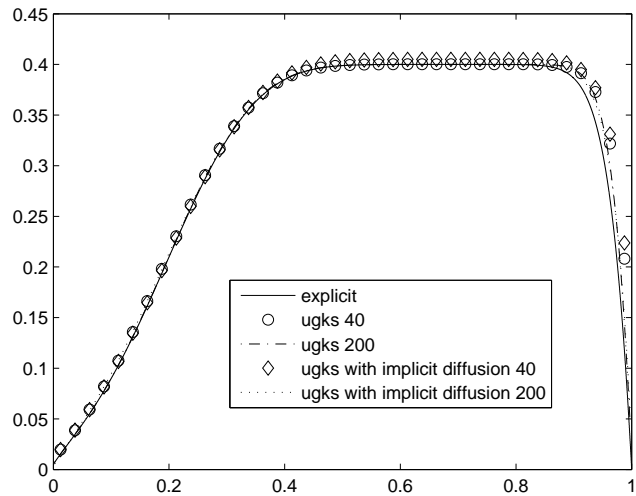}}
\centerline{\includegraphics[width=0.6\textwidth]{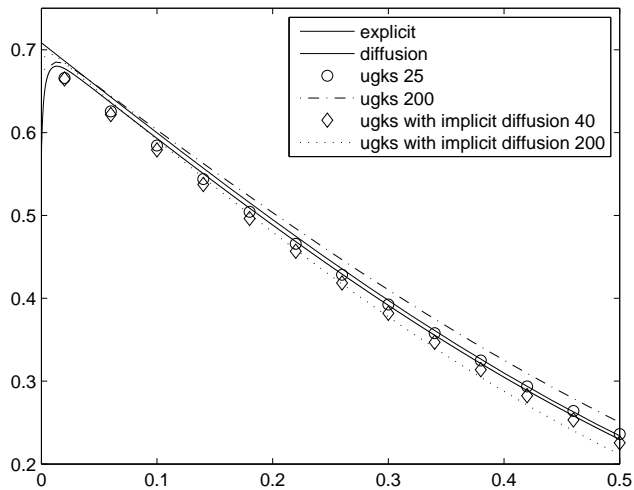}}
\caption{\label{fig:diff_imp_JPT2_E4}Comparison of the UGKS with/without
  implicit diffusion with a reference solution for long time cases. Data of of examples
  3 (Intermediate regime with a variable scattering frequency and a source term, top)
  and 5 (intermediate regime with non-isotropic boundary
conditions, bottom).}
\end{figure}

In figure~\ref{fig:diff_imp_JPT2_E4}, we compare our schemes for long
time cases.  The data are those of examples 3 (see
section~\ref{subsec:isotropic}), and 5 (see
section~\ref{subsec:res_bc}). We observe that for the coarse mesh, both schemes are very
close. For example 3, the time step of UGKS-ID is 25 times 
as large as for UGKS-ED, while it is 16 times as large for example 5.
For the fine mesh, UGKS-ID
requires $cfl=0.1$ (it is unstable for larger CFL), and the time step
is here 100 times as large as for UGKS-ED in example 3, and 10 times
as large as for example 5. For example 5, UGKS-ID is different from the
reference solution, but this difference is as large as the difference
observed for UGKS-ED.

This study shows that the UGKS-ID can be efficiently used when
the parabolic CFL condition is too much restrictive, with the same
accuracy as the non-modified UGKS-ED.

%--------------------------------------------------------------------------
\section{Conclusion}
%--------------------------------------------------------------------------

In this paper, we have shown that the unified gas kinetic scheme,
originally designed for rarefied gas dynamics problems, can be applied
to other kinetic equations, like radiative transfer models. Moreover,
the UGKS has been shown to be asymptotic preserving in the diffusion
limit of such equations, as well as in the free transport
regime. This scheme turns out to be an efficient multiscale method for
kinetic problems, with a wide range of applications. In addition, we
have shown that the UGKS can be simply modified to account for boundary
layers and to obtain an implicit scheme in the diffusion
limit. Finally, we have suggested an extension of the UGKS to general
collision operator (in a non relaxation form) like that of neutron
transport, that is still to be tested. This extension preserves the AP property of the method
under some conditions. 

It would be interesting to rigorously prove the stability of the UGKS, in
particular to determine an explicit CFL condition: we believe that the
energy method we proposed for another AP scheme in~\cite{LM2010} could
be applied. For practical applications, we would like to extend this
work to multidimensional problems. This should not be difficult, since
it has already been done in rarefied gas dynamics by Xu and Huang in~\cite{HXY2012}.

\bibliographystyle{plain}

\bibliography{biblio}

\begin{thebibliography}{10}

\bibitem{BLM2008}
M.~Bennoune, M.~Lemou, and L.~Mieussens.
\newblock Uniformly stable numerical schemes for the {B}oltzmann equation
  preserving the compressible {N}avier-{S}tokes asymptotics.
\newblock {\em J. Comput. Phys.}, 227(8):3781--3803, 2008.

\bibitem{BCLM}
C.~Buet, S.~Cordier, B.~Lucquin-Desreux, and S.~Mancini.
\newblock Diffusion limit of the {L}orentz model: asymptotic preserving
  schemes.
\newblock {\em M2AN Math. Model. Numer. Anal.}, 36(4):631--655, 2002.

\bibitem{BDF2012}
C.~Buet, B.~Despr\'es, and E.~Franck.
\newblock Design of asymptotic preserving finite volume schemes for the
  hyperbolic heat equation on unstructured meshes.
\newblock {\em Numerische Mathematik}, 122:227--278, 2012.

\bibitem{CGL2008}
J.~A. Carrillo, T.~Goudon, and P.~Lafitte.
\newblock Simulation of fluid and particles flows: asymptotic preserving
  schemes for bubbling and flowing regimes.
\newblock {\em J. Comput. Phys.}, 227(16):7929--7951, 2008.

\bibitem{CGLV2008}
J.~A. Carrillo., T.~Goudon, P.~Lafitte, and F.~Vecil.
\newblock Numerical schemes of diffusion asymptotics and moment closures for
  kinetic equations.
\newblock {\em J. Sci. Comput.}, 36(1):113--149, 2008.

\bibitem{FJ2010}
F.~Filbet and S.~Jin.
\newblock {A class of asymptotic-preserving schemes for kinetic equations and
  related problems with stiff sources}.
\newblock {\em Journal of Computational Physics}, 229(20):7625--7648, 2010.

\bibitem{Gosse_2011}
L.~Gosse.
\newblock {Transient radiative tranfer in the grey case: well-balanced and
  asymptotic-preserving schemes built on Case's elementary solutions}.
\newblock {\em Journal of Quantitative Spectroscopy and Radiative Transfer},
  112:1995--2012, 2011.

\bibitem{GT2002}
L.~Gosse and G.~Toscani.
\newblock An asymptotic-preserving well-balanced scheme for the hyperbolic heat
  equations.
\newblock {\em C. R. Math. Acad. Sci. Paris}, 334(4):337--342, 2002.

\bibitem{GT2003}
L.~Gosse and G.~Toscani.
\newblock Space localization and well-balanced schemes for discrete kinetic
  models in diffusive regimes.
\newblock {\em SIAM J. Numer. Anal.}, 41(2):641--658, 2003.

\bibitem{HXY2012}
J.C. Huang, K.~Xu, and P.B. Yu.
\newblock A unified gas-kinetic scheme for continuum and rarefied flows ii:
  Multi-dimensional cases.
\newblock {\em Communications in Computational Physics}, 3(3):662--690, 2012.

\bibitem{jin_sisc_1999}
S.~Jin.
\newblock Efficient asymptotic-preserving ({AP}) schemes for some multiscale
  kinetic equations.
\newblock {\em SIAM J. Sci. Comput.}, 21(2):441--454, 1999.

\bibitem{JL1991}
S.~Jin and C.~D. Levermore.
\newblock The discrete-ordinate method in diffusive regimes.
\newblock {\em Transport Theory Statist. Phys.}, 20(5-6):413--439, 1991.

\bibitem{JL1993}
S.~Jin and C.~D. Levermore.
\newblock Fully discrete numerical transfer in diffusive regimes.
\newblock {\em Transport Theory Statist. Phys.}, 22(6):739--791, 1993.

\bibitem{jin_levermore}
S.~Jin and C.~D. Levermore.
\newblock Numerical schemes for hyperbolic conservation laws with stiff
  relaxation terms.
\newblock {\em J. Comput. Phys.}, 126:449, 1996.

\bibitem{JP2000}
S.~Jin and L.~Pareschi.
\newblock Discretization of the multiscale semiconductor {B}oltzmann equation
  by diffusive relaxation schemes.
\newblock {\em J. Comput. Phys.}, 161(1):312--330, 2000.

\bibitem{JP2000b}
S.~Jin and L.~Pareschi.
\newblock Asymptotic-preserving ({AP}) schemes for multiscale kinetic
  equations: a unified approach.
\newblock In {\em Hyperbolic problems: theory, numerics, applications, Vol. I,
  II (Magdeburg, 2000)}, volume 141 of {\em Internat. Ser. Numer. Math., 140},
  pages 573--582. Birkh\"auser, Basel, 2001.

\bibitem{JPT1998}
S.~Jin, L.~Pareschi, and G.~Toscani.
\newblock Diffusive relaxation schemes for multiscale discrete-velocity kinetic
  equations.
\newblock {\em SIAM J. Numer. Anal.}, 35(6):2405--2439, 1998.

\bibitem{JPT2000}
S.~Jin, L.~Pareschi, and G.~Toscani.
\newblock Uniformly accurate diffusive relaxation schemes for multiscale
  transport equations.
\newblock {\em SIAM J. Numer. Anal.}, 38(3):913--936, 2000.

\bibitem{klar_sinum_1998}
A.~Klar.
\newblock An asymptotic-induced scheme for nonstationary transport equations in
  the diffusive limit.
\newblock {\em SIAM J. Numer. Anal.}, 35(3):1073--1094, 1998.

\bibitem{klar_sinum_1999}
A.~Klar.
\newblock An asymptotic preserving numerical scheme for kinetic equations in
  the low {M}ach number limit.
\newblock {\em SIAM J. Numer. Anal.}, 36(5):1507--1527, 1999.

\bibitem{klar_sisc_1999}
A.~Klar.
\newblock A numerical method for kinetic semiconductor equations in the
  drift-diffusion limit.
\newblock {\em SIAM J. Sci. Comput.}, 20(5):1696--1712, 1999.

\bibitem{KS2001}
A.~Klar and C.~Schmeiser.
\newblock Numerical passage from radiative heat transfer to nonlinear diffusion
  models.
\newblock {\em Math. Models Methods Appl. Sci.}, 11(5):749--767, 2001.

\bibitem{KU2002}
A.~Klar and A.~Unterreiter.
\newblock Uniform stability of a finite difference scheme for transport
  equations in diffusive regimes.
\newblock {\em SIAM J. Numer. Anal.}, 40(3):891--913, 2002.

\bibitem{LS2012}
P.~Lafitte and G.~Samaey.
\newblock Asymptotic-preserving projective integration schemes for kinetic
  equations in the diffusion limit.
\newblock {\em SIAM Journal on Scientific Computing}, 34(2):A579--A602, 2012.

\bibitem{LM1989}
A.~W. Larsen and J.~E. Morel.
\newblock Asymptotic solutions of numerical transport problems in optically
  thick, diffusive regimes. {II}.
\newblock {\em J. Comput. Phys.}, 83(1):212--236, 1989.

\bibitem{LMM1987}
A.~W. Larsen, J.~E. Morel, and W.~F.~Miller Jr.
\newblock Asymptotic solutions of numerical transport problems in optically
  thick, diffusive regimes.
\newblock {\em J. Comput. Phys.}, 69(2):283--324, 1987.

\bibitem{Lemou_relaxedAP_2010}
M.~Lemou.
\newblock {Relaxed micro-macro schemes for kinetic equations}.
\newblock {\em C.R Acad. Sci.}, Serie I, 348(7-8):455--460, 2010.

\bibitem{LM2013}
M.~Lemou and F.~M\'ehats.
\newblock Micro-macro schemes for kinetic equations including boundary layers.
\newblock {\em SIAM Journal on Numerical Analysis}.
\newblock to appear.

\bibitem{LM2007}
M.~Lemou and L.~Mieussens.
\newblock A new asymptotic preserving scheme based on micro-macro formulation
  for linear kinetic equations in the diffusion limit.
\newblock {\em SIAM J. Sci. Comput.}, 31(1):334--368, 2008.

\bibitem{LM2010}
J.~Liu and L.~Mieussens.
\newblock Analysis of an asymptotic preserving scheme for linear kinetic
  equations in the diffusion limit.
\newblock {\em SIAM Journal on Numerical Analysis}, 48(4):1474--1491, 2010.

\bibitem{MELD_2008}
R.~G. McClarren, T.~M. Evans, R.~B. Lowrie, and J.~D. Densmore.
\newblock Semi-implicit time integration for pn thermal radiative transfer.
\newblock {\em Journal of Computational Physics}, 227:7561–7586, 2008.

\bibitem{NP1998}
G.~Naldi and L.~Pareschi.
\newblock Numerical schemes for kinetic equations in diffusive regimes.
\newblock {\em Appl. Math. Lett.}, 11(2):29--35, 1998.

\bibitem{PP2007}
Sandra Pieraccini and Gabriella Puppo.
\newblock Implicit–explicit schemes for bgk kinetic equations.
\newblock {\em Journal of Scientific Computing}, 32:1--28, 2007.

\bibitem{TX2006}
Manuel Torrilhon and Kun Xu.
\newblock Stability and consistency of kinetic upwinding for
  advection‚Äìdiffusion equations.
\newblock {\em IMA Journal of Numerical Analysis}, 26(4):686--722, October
  2006.

\bibitem{wesseling}
P.~Wesseling.
\newblock {\em Principles of Computational Fluid Dynamics}, volume~29 of {\em
  Series: Springer Series in Computational Mathematics}.
\newblock Springer, 2001.

\bibitem{Xu_2001}
K.~Xu.
\newblock A gas-kinetic {BGK} scheme for the navier–stokes equations and its
  connection with artificial dissipation and godunov method.
\newblock {\em Journal of Computational Physics}, 171(1):289 -- 335, 2001.

\bibitem{XH2010}
K.~Xu and J.-C. Huang.
\newblock A unified gas-kinetic scheme for continuum and rarefied flows.
\newblock {\em J. Comput. Phys.}, 229:7747--7764, 2010.

\end{thebibliography}

%----- FIGURES ----------------
\newpage

%----- various regimes, isotropic BC
\clearpage

\clearpage

\clearpage

\clearpage

%----- Boundary layer
\clearpage

\clearpage

\clearpage

%----- Implicit diffusion
\clearpage

\clearpage

\end{document}